\documentclass[11pt]{article}
\usepackage{amscd,amssymb,stmaryrd}

\usepackage{graphicx}
\usepackage{subcaption}
\usepackage[utf8]{inputenc}
\usepackage[export]{adjustbox}
\usepackage{wrapfig}

\usepackage{amsthm}
\usepackage{bbm}
\usepackage{cite}
\usepackage{color}
\usepackage{mathtools}
\usepackage{diagbox} 
\usepackage{hyperref}
\usepackage[english]{babel}
\usepackage{booktabs}

\usepackage{float}
\usepackage[linesnumbered,ruled,vlined, ruled,vlined]{algorithm2e}

\theoremstyle{definition}

\usepackage{tikz}
\usepackage{pifont}%

\usepackage{diagbox}
\usetikzlibrary{decorations.pathreplacing,calligraphy}
\usepackage{pgfplots}
\pgfplotsset{compat=1.11}
\usepgfplotslibrary{groupplots}

\usetikzlibrary{shapes.misc}

\tikzset{cross/.style={cross out, draw=black, minimum size=2*(#1-\pgflinewidth), inner sep=0pt, outer sep=0pt},
cross/.default={1pt}}

\usepackage[linesnumbered,ruled,vlined, ruled,vlined]{algorithm2e}
\usepackage{amsmath}
\usepackage{lineno}
\usepackage{subcaption}
\usepackage{tikz}
\usepackage{xcolor}
\usepgfplotslibrary{fillbetween}
\usepackage{float}
\usepackage{pgfplots}
\pgfplotsset{compat=newest} 

\title{Finite Element Representation Network (FERN) for Operator Learning with a Localized Trainable Basis}
\author{Zecheng Zhang \footnote{Department of Applied and Computational Mathematics and Statistics, University of Notre Dame, Notre Dame, IN 46556, USA. (Email: zecheng.zhang.math@gmail.com)},
Hao Liu \footnote{Department of Mathematics, Hong Kong Baptist University, Hong Kong, China. (Email: haoliu@hkbu.edu.hk)},
Guosheng Fu\footnote{Department of Applied and Computational Mathematics and Statistics, University of Notre Dame, Notre Dame, IN 46556, USA. (Email: gfu@nd.edu)},
Hayden Schaeffer \footnote{Department of Mathematics, UCLA, Los Angeles, CA 90095, USA. (Email: hayden@math.ucla.edu)},
Guang Lin \footnote{Department of Mathematics and Mechanical Engineering, Purdue University, West Lafayette, IN 47907, USA. (Email: guanglin@purdue.edu)}.
}
\date{}

\begin{document}
\maketitle
 
\begin{abstract}
We propose a finite-element local basis-based operator learning framework for solving partial differential equations (PDEs). 
Operator learning aims to approximate mappings from input functions to output functions, where the latter are typically represented using basis functions. 
While non-learnable bases reduce training costs, learnable bases offer greater flexibility but often require deep network architectures with a large number of trainable parameters. 
Existing approaches typically rely on deep global bases; however, many PDE solutions exhibit local behaviors such as shocks, sharp gradients, etc., and in parametrized PDE settings, these localized features may appear in different regions of the domain across different training and testing samples. 
Motivated by the use of local bases in finite element methods (FEM) for function approximation, we develop a shallow neural network architecture that constructs adaptive FEM bases. 
By adopting suitable activation functions, such as ReLU, the FEM bases can be assembled exactly within the network, introducing no additional approximation error in the basis construction process. 
This design enables the learning procedure to naturally mimic the adaptive refinement mechanism of FEM, allowing the network to discover basis functions tailored to intrinsic solution features such as shocks. 
The proposed learnable adaptive bases are then employed to represent the solution (output function) of the PDE. This framework reduces the number of trainable parameters while maintaining high approximation accuracy, effectively combining the adaptivity of FEM with the expressive power of operator learning. 
To evaluate performance, we validate the proposed method on seven families of PDEs with diverse characteristics, demonstrating its accuracy, efficiency, and robustness.

\end{abstract}

\section{Introduction}
Deep neural operators (DNOs) are a class of operator learning (OL) approaches \cite{lu2021learning,lu2022comprehensive,zhang2023belnet,lin2023b,zhang2024d2no, lin2021operator,howard2022multifidelity, meuris2023machine,li2020fourier,li2022fourier,wen2022u,lanthaler2024discretization,kovachki2021universal,goswami2022physics,hao2024newton, de2022cost} that aims to approximate an operator $G$, which maps an input function $u$ to an output function $v = G(u)$, by utilizing deep neural networks architectures. 
Since many problems in scientific computing can be formulated as operator approximations, OL has emerged as a powerful alternative technique within scientific machine learning (SciML) \cite{raissi2019physics, lu2021deepxde, lu2021learning, schaeffer2017sparse, zhang2019convergence, schaeffer2017sparse, schaeffer2013sparse, brunton2016discovering} for solving computational problems. One of their advantages over classical techniques is their faster inference time and their ability to adapt to problem-specific tasks. It has been shown theoretically that deep neural network are adaptive to data structures and model smoothness \cite{liu2021besov,nakada2020adaptive,chen2022nonparametric,suzuki2018adaptivity,havrilla2024understanding}.
Recent advances in multi-operator learning (MOL) \cite{liu2023prose,sun2024towards,yang2023prompting,yang2023context,yang2025fine,subramanian2024towards,negrini2025multimodal, yu2024nonlocal,ye2024pdeformer,liu2025bcat, jollie2025time, cao2024vicon, liu2024prosefd},  which employ a single network architecture to approximate multiple operators simultaneously, have significantly expanded the scope of OL. 
This development enables the application of OL to large-scale engineering problems \cite{pathak2022fourcastnet, guibas2021adaptive, zhu2023fourier, deng2022openfwi, kurth2023fourcastnet, mao2024disk2planet, wang2025learning, bodnar2025foundation}, positioning it as one of the few SciML frameworks already adopted in practical applications.

The study of operator learning started in the early 1990s with the works \cite{chen1995universal, chen1993approximations}, which established a universal approximation theorem for operators and laid the theoretical foundation of the field.  In recent years, a variety of DNOs framework have been developed \cite{lu2021learning,zhang2023belnet,lin2023b,zhang2024d2no,zhang2024modno, lin2021operator,howard2022multifidelity, meuris2023machine, goswami2022physics,kumar2025synergistic, li2020fourier,li2023fourier,li2022fourier,wen2022u}.
For example, DeepONet is the first modern deep neural operator and was able to show applicability to spatiotemporal systems. However, the vanilla version of DeepONet is not invariant to the discretization of the input function; in particular, the input function cannot be discretized on different meshes. This limitation was partially resolved in \cite{zhang2023belnet} under suitable assumptions. Nevertheless, the discretization size was still required to remain the same across different functions. This restriction was improved on both numerically and theoretically in \cite{zhang2024d2no, zhang2024modno}, through the introduction of distributed learning strategies.
Significant theoretical progress \cite{liu2024deep, lanthaler2023parametric, bhattacharya2021model, lanthaler2023operator, liu2024neural, zhang2025coefficient, lanthaler2023parametric, kovachki2021universal, lanthaler2022error, liu2025generalization} has been made to analyze properties such as convergence, generalization, and expressivity of neural operators, further solidifying OL as an important topic in SciML.

A central mathematical motivation for the design of DNOs is to represent the output function $v$ as a linear combination of learnable basis functions (represent the basis functions using neural networks) with corresponding learnable coefficients. 
This idea is commonly used in numerical methods, such as the finite element method and spectral methods, which express functions as expansions over basis functions. Within OL, this idea was first illustrated in the shallow neural network universal approximation results of \cite{chen1993approximations, chen1995universal}, and was later rigorously established in \cite{liu2024neural}, which provided complete neural scaling laws, including error convergence with respect to network size and generalization error with respect to the number of training samples.  

An example of such a structure is the DeepONet \cite{lu2021learning, lin2023b}, which employs a trunk network to learn basis functions and a branch network to learn the coefficients. 
Another widely used framework is the FNO \cite{li2020fourier,li2023fourier,li2022fourier,wen2022u} which adopts the iterative structure $u_{n+1} = \mathcal{F}^{-1}(\mathcal{K}_n(\mathcal{F}(u_n))) + W u_n,$ 
where $\mathcal{F}$ denotes the Fourier transform, $W$ is a learnable weight matrix, and $\mathcal{K}_n$ is a learned convolution. 
In this formulation, although part of the basis is fixed as Fourier and inverse Fourier modes, the term $W u$ still represents a linear combination of learned column basis functions weighted by coefficients determined by $u$.   More recent architectures that leverage attention mechanisms \cite{liu2023prose, liu2025bcat, subramanian2024towards, efendiev2022efficient} to construct DNOs can also be interpreted within a similar framework, as attention effectively implements a data-dependent basis expansion with learned coefficients.

Both numerical evidence \cite{lu2021learning, lin2023b} and theoretical analysis \cite{liu2024neural} demonstrate that, in order to achieve strong approximation performance, the networks used to learn the basis functions are typically deep, which results in globally supported basis functions, i.e., their support spans the entire domain of the input function $u$. Such approaches can be viewed as spectral methods \cite{meuris2023machine}, but with the important distinction that the basis functions are learnable and non-polynomial (depending on the choice of the neural network activation function).
While global basis functions offer advantages such as smoothness and high accuracy for problems with globally coherent structures, they are less effective for problems involving localized phenomena, sharp gradients, or discontinuities. This is common in many PDEs with shocks, interfaces, or boundary layers. 
In such cases, locally supported basis functions, e.g., those used in the finite element or finite volume methods, generally yield better performance by more accurately resolving localized features.
Moreover, OL aims to predict the solutions of a parametrized family of PDEs, implying that the localized features of different samples are often evenly distributed across the solution landscapes. Consequently, local bases are sufficient to capture these features effectively, eliminating the need for global bases.

By contrast, the finite element method (FEM), originally introduced in \cite{hrennikoff1941solution, feng1965finite, courant1994variational}, provides a rigorous framework for constructing locally supported basis functions \cite{Ciarlet2002Finite, Johnson1987, BrennerScott2008} utilizing in solving PDE. FEM is one of the main approach employed in numerical methods since it is capable of handling a wide range of PDE with complex geometries through the use of unstructured meshes \cite{Zienkiewicz13, Bathe1996}. 
In the most common setting, the basis elements are piecewise linear Lagrange basis functions, often referred to as hat functions, defined on elements of a triangulation of the computational domain. These locally supported basis functions enable accurate resolution of localized features.

The ReLU activation function \cite{yarotsky2017error, liu2024neural} is capable of constructing piecewise linear functions, which can be used for interpolation and for defining basis functions with local support. 
In particular, with only two degrees of freedom (trainable parameters), one can construct a one-dimensional hat function. Motivated by this observation, we propose to 2-parameter ReLU networks to construct local FEM-type basis functions for operator learning. 
Since each basis is parameterized by trainable variables, the approach can be regarded as an adaptive FEM method. We validate the proposed framework on a broad family of PDEs that are traditionally solved by FEM, and demonstrate that the newly designed neural operators achieve comparable or even improved accuracy while requiring significantly lower computational cost.
Since our approach is directly inspired by the finite element representation of functions, we refer to it as the Finite Element Representation Network (FERN).

We summarize our contributions as follows:
\begin{enumerate}
\item We propose a novel operator learning framework based on finite-element local bases, where adaptive FEM bases are constructed directly through a shallow neural network. 
The use of local bases is crucial, as PDE solutions may exhibit localized behaviors such as shocks, sharp gradients, or fast decays. 
Moreover, in parametrized PDE problems, these local features may appear at different spatial locations across samples, making local bases essential for accurately capturing such variations while maintaining efficiency and interpretability.

\item By adopting suitable activation functions such as ReLU, the FEM bases can be assembled exactly within the network, introducing no approximation error in the basis construction process. 
The learning procedure naturally mimics the adaptive refinement mechanism of adaptive FEM, enabling the network to discover basis functions tailored to the underlying solution structures.  
\item The proposed framework achieves high prediction accuracy with a reduced number of trainable parameters, owing to the low cost of constructing the bases. This provides an effective approach to employing shallow network structures while maintaining strong predictive performance. 
Extensive experiments on PDEs with distinct characteristics, including solutions whose landscapes exhibit bumps, demonstrate the method’s accuracy, efficiency, and robustness.
\end{enumerate}
The rest of the manuscript is organized as follows.
We introduce our methodology in Section \ref{sec_method_review}, and will verify our methods and present the numerical evidence in Section \ref{sec_numerilca}.
We finally discuss some future work in Section \ref{sec_conclusion}.

\section{Review and methodology}
\label{sec_method_review}
\subsection{Operator learning}
Let us denote the operator $G: U \rightarrow V$ as the target operator, where $U$ and $V$ are function spaces defined on domains $\Omega_u$ and $\Omega_v$ with dimensions $d_1$ and $d_2$, respectively. 
The objective of operator learning is to design a neural network structure $G_{\theta}$, with $\theta$ representing all trainable parameters, to approximate $G$. More precisely, we aim to achieve
\begin{align*}
    G(u)(x) := v(x) \approx G_{\theta}(\hat{u})(x),  
\end{align*}
for $u \in U$ and $x \in \Omega_v$, where $v = G(u)$ and $\hat{u}$ is a discretization of $u$.  

A useful construction of $G_{\theta}$ is inspired by the classical numerical analysis, which first approximates $v \in V$ as
\begin{align*}
    v(x) = \sum_{k=1}^N c_k \, \phi_k(x),
\end{align*}
where $c_k$ are coefficients and $\phi_k$ are basis functions. 
In standard numerical methods, the basis functions $\phi_k$ are chosen independently of the solution $v(x)$ and are fixed \emph{a priori}, for example, as finite element basis functions. 
As a concrete example, in the low-order continuous finite element method, the function space $V$ is defined as the space of continuous, piecewise linear functions over a fixed computational mesh. Each basis function $\phi_k$ is a {hat function} associated with node $k$ of the mesh, taking the value one at node $k$ and zero at all other nodes. 
This construction represents the function $v$ as a linear combination of locally supported basis functions, yielding a flexible yet structured approximation framework.

In operator learning, the coefficients $c_k$ must encode both $u$ and $G$. Consequently, one can approximate the operator as
\begin{align}
    G(u)(x) \approx \sum_{k=1}^N c_k(u)\,\phi_k(x),
    \label{eqn_operator_approx}
\end{align}
where $c_k: U \rightarrow \mathbb{R}$ are functionals on $U$ encoding both $u$ and $G$. To realize this approximation, one designs network structures $\hat{c}_k$ and $\hat{\phi}_k$ as surrogates for $c_k$ and $\phi_k$, respectively.  
In the remainder of this work, we adopt the terminology of \cite{lu2021learning}, referring to the networks $\hat{c}_k$ as branch networks and $\hat{\phi}_k$ as trunk networks.

This framework was first studied mathematically in \cite{chen1995universal}, whose authors analyzed shallow approximations of $\hat{c}_k$ and $\hat{\phi}_k$, although the convergence rates and parameter estimates for encoding $u$, encoding the operator, and constructing the basis $\hat{\phi}_k$ remained unclear. Subsequently, \cite{lu2021learning, lu2022comprehensive, meuris2023machine} proposed deep computational extensions of this shallow structure, significantly improving its empirical performance but without providing theoretical guarantees. 
Finally, the work of \cite{liu2024neural} rigorously established convergence and generalization results for this operator learning framework, which was the first explicit estimates of the parameter complexity required for $u$ encoding, operator encoding, and basis construction. 

Notably, \cite{liu2024neural} shows that the depth of the network $\hat{\phi}_k$ required to achieve an error of at most $\varepsilon > 0$ in the $L^{\infty}$ norm scales as $-\log(\varepsilon)$. 
This implies that a deep architecture is necessary, which can be computationally expensive. 
This theoretical result is consistent with numerical observations: shallow structures for $\hat{\phi}_k$ often fail to yield satisfactory computational performance (see Section \ref{sec_numerilca} for examples).
A central challenge is to reduce the cost of constructing the basis while maintaining high prediction accuracy.

Additionally, the suggested total parameters count \cite{liu2024neural} for $\hat{c}_k$ is of order $(\sqrt{c_U}\varepsilon^{-(d_2+1)c_U})(c_U^2\log c_U+c_U^2\log(\varepsilon^{-1}))$ with $C_u = O(\varepsilon^{-d_1})$.
This estimate dominates the parameter count for trunk network $\hat{\phi}_k$, which is of order $-\log(\varepsilon)$. 
However, in numerical experiments, the parameter count of the branch network typically dominates that of the trunk network to ensure accurate performance.
This raises an important question: how can we design a structure that achieves high predictive accuracy while requiring fewer trainable parameters in $\hat{\phi}_k$, so that the network structures $\hat{c}_k$ and $\hat{\phi}_k$ are consistent with theoretical expectations.

\subsection{Finite element methods}
The finite element method is a powerful numerical technique for solving partial differential equations (PDEs). One of its key strengths lies in its ability to operate on unstructured meshes that can accurately represent complex geometries. In this framework, the computational domain $\Omega$ 
is discretized into a triangulation $\mathcal{T}_h=\{K\}$, where each 
$K$ denotes an individual mesh element. 
Based on the continuity requirements of the problem, local finite element spaces are then constructed over this mesh. For example, the classical continuous piecewise linear Lagrange finite element space on a simplicial mesh is defined as
\begin{align}
\label{p1}
V_h(\mathcal{T}_h) = \{ v_h \in C^0(\Omega) : v_h|_K \in P_1(K) \text{ for all } K \in \mathcal{T}_h \},    
\end{align}
where $P_1(K)$ represents the space of linear polynomials over the simplicial element $K$.

In particular, adaptive finite element methods (AFEM) dynamically refine the mesh to better resolve local features of the solution according to various adaptivity criteria~\cite{Mark00,Bangerth03,Nochetto24}. 
For instance, let $\Omega \subset \mathbb{R}^2$ be a polygonal domain. 
If $u \in H^{s}(\Omega)$ with $0 < s < 2$, then there exists an admissible triangulation 
$\mathcal{T}_h$ with $\mathcal{O}(N)$ vertices such that~\cite{Binev02}
\begin{equation}\label{eq:best_approximation}
    \inf_{v_h \in V_h(\mathcal{T}_h)} 
    \|u - v_h\|_{L^2(\Omega)} 
    \le C\, N^{-s/2}\, |u|_{H^{s}(\Omega)},
\end{equation}
where $V_h(\mathcal{T}_h)$ denotes the space of continuous, piecewise linear 
finite element functions defined in~\eqref{p1}. 
The constant $C>0$ depends only on the domain $\Omega$ and on the shape-regularity of the mesh family.

The remarkable performance of adaptive FEM motivates the idea of constructing finite element partitions and basis functions directly through neural networks, where the adaptive process of finding optimal basis functions is embedded and optimized within the network training. To this end, we utilize the ReLU activation function, which serves as a fundamental building block in designing such FEM-inspired neural architectures.

\subsection{ReLU functions}
The ReLU function, defined as $\sigma(x) = \max\{x,0\}$, is a widely used activation function for introducing nonlinearity in neural networks. A notable feature of ReLU is its ability to generate piecewise linear functions. Consider the interval $[0,1]$ partitioned by a mesh $x_0 = 0 < x_1 < \cdots < x_n = 1$. Given prescribed slopes $k_i$ on each subinterval $[x_{i-1},x_i]$ for $i=1,\dots,n$, the corresponding piecewise linear function can be expressed as  
\begin{align}
    f(x) = k_1 \sigma(x) + \sum_{j=1}^{n-1} (k_{j+1}-k_j)\,\sigma(x-x_j).
    \label{eqn_pwl_relu}
\end{align}
This representation highlights the flexibility of ReLU in constructing piecewise linear functions and serves as the motivation for employing ReLU to build finite element basis functions.
In fact, this connection was rigorously studied in \cite{HeLiXuZheng2020}, showing that continuous piecewise linear (CPWL) finite element functions can be represented exactly by deep ReLU neural networks. Moreover, they establish lower bounds on the depth required (e.g.\ at least two hidden layers for 
$d\ge 2$ with $d$ denoting the spatial dimension) and discuss how this equivalence underlies the expressive power of ReLU networks in approximating classical finite element spaces.

\subsection{FEM basis construction}
\begin{figure}[H]
    \centering
    \includegraphics[scale = 0.4]{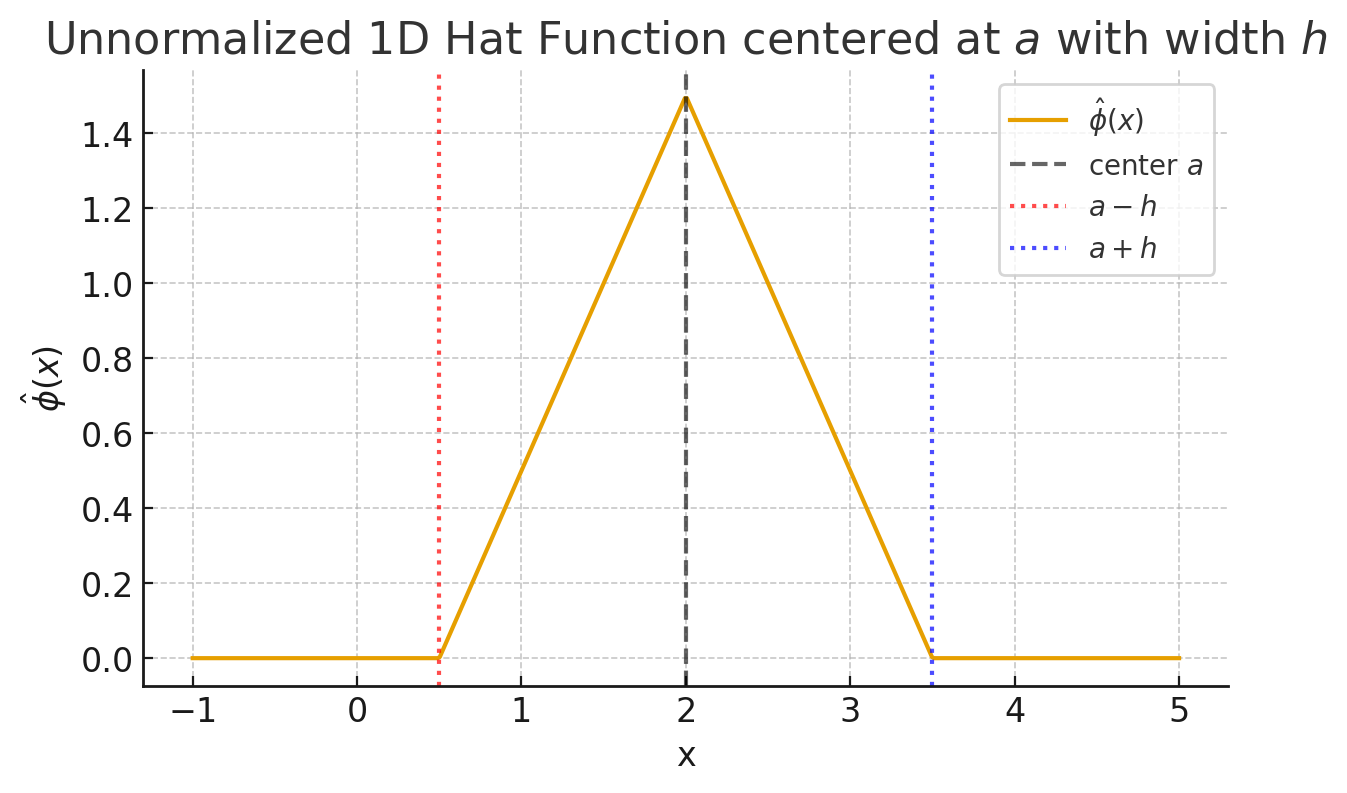}
    \caption{A hat function demonstration.}
    \label{fig_hat_function}
\end{figure}

We construct the network following the structure in Equation~(\ref{eqn_operator_approx}). 
Instead of employing a deep architecture (which typically requires a large number of trainable parameters) to build the basis functions $\phi_k(x)$, we consider constructing the basis $\phi_k$ using FEM hat functions. 
Specifically, a hat function (illustrated in Figure~\ref{fig_hat_function}) is defined as
\begin{align*}
    p_{a, h}(x) =
\begin{cases}
x - (a-h), & a-h < x \leq a, \\
-x + (a+h), & a < x < a+h, \\
0, & \text{otherwise}.
\end{cases}
\end{align*}
By Equation (\ref{eqn_pwl_relu}), this hat function $p_{a, h}$ can be represented exactly by the ReLU functions as
\begin{align*}
    p_{a, h}(x) = \sigma\left(x-(a-h)) - 2\sigma(x - a) +\sigma(x-(a+h)\right).
\end{align*}
Notably, to represent each hat function, only 2 parameters are used.
We then formulate our approximation as
\begin{align*}
    G_{\theta}(u)(x) = \sum_{k = 1}^N \hat{c}_{k, \theta_k}(\hat{u} ) p_{a_k, h_k}(x) = \sum_{k = 1}^N \hat{c}_{k, \theta_k} \left(\sigma(x-(a_k-h_k)) - 2\sigma(x - a_k) +\sigma(x - (a_k+h_k) ) \right),
\end{align*}
where $\hat{c}_{k, \theta_k}$ denotes the neural network architecture with parameters $\theta_k$, $\theta$ denotes the collection of all $\theta_k$, and the trainable parameters $a_k, h_k$ used to construct the hat function basis.

As a result, the FEM basis can be constructed exactly as a combination of ReLU functions, introducing no approximation error in the basis construction process. During network training, the centers $a_k$ and widths $h_k$ are treated as learnable parameters, effectively mimicking the adaptive refinement process in FEM and enabling the network to accurately capture the output functions (e.g., the PDE solutions). 
Examples illustrating this behavior can be found in Sections~\ref{sec_fk_example}, \ref{sec_aggregation}, \ref{sec_keller}, and~\ref{sec_kdv}, where the proposed method successfully captures bumps or sharp peaks (singular features) in the solution landscapes, demonstrating its capability to represent non-smooth behaviors effectively.

\section{Numerical examples}
\label{sec_numerilca}
In this section, we present a series of numerical examples to demonstrate the performance of the proposed method. We consider seven PDEs that exhibit a wide range of physical and mathematical properties, including gradient-flow dynamics, dispersive behavior, and nonlinear advection–diffusion. Specifically, the following PDEs are studied: the Allen–Cahn (AC) equation, the Cahn–Hilliard (CH) equation, the Fokker–Planck (FP) equation, the Aggregation–Diffusion (AD) equation, the Keller–Segel (KS) equation, the Korteweg–de Vries (KdV) equation, and the viscous Burgers equation. Among these, the Allen–Cahn, Cahn–Hilliard, Fokker–Planck, Aggregation–Diffusion, and Keller–Segel equation are examples of gradient-flow systems that evolve toward equilibrium under an energy-dissipation principle, whereas the KdV equation represents a dispersive system, and the viscous Burgers equation is a dissipative nonlinear advection–diffusion model.
To show the adaptivity of the methods and later shows that our methods mimic the adaptive FEM, in FK, AD, KS and KdV example, we created either bumps, shocks, or fast decay in the solution landscape, which is capture by our method. 
Furthermore, to demonstrate the adaptivity of the proposed method—and to show that it effectively mimics the adaptive FEM, which explains its high numerical accuracy and efficiency—we design experiments in the FK, AD, KS, and KdV examples where the solution landscapes exhibit bumps, shocks, or rapid decays. 

To fairly evaluate the proposed method, we compare our method against the standard DeepONet, which can be viewed as using non-polynomial global bases~\cite{meuris2023machine}. 
For a fair comparison, the coefficient branch network $\hat{c}_{k,\theta_k}(\cdot)$ structure is kept identical between our method and the standard DeepONet. 
However, the DeepONet employs a deep (and, for additional comparison with our shallow FEM network, a shallow) fully connected architecture for its trunk network, while our method uses the proposed shallow structures to construct the FEM local polynomial. 
The number of basis functions $N$ is also kept the same across all comparison tests. In addition, we include comparisons with the POD-based approach, where the POD modes of the output functions are used as the trunk network in examples where such a method is applicable.

\subsection{Summary of Results} 
\textbf{Accuracy and efficiency.}
Owing to the use of hat functions in constructing the basis, the proposed FEM method FERN substantially reduces the number of training parameters. 
Despite having significantly fewer parameters than DeepONet with an identical structure (i.e., the same number of basis functions $N$ and the same coefficient networks $c_{k,\theta_k}$), we observe comparable accuracy across all examples, and in $10$ out of the $13$ comparison tests, the FERN achieves improved accuracy with smaller prediction variance across all testing samples. 
In contrast, when a two-layer fully connected network is used as the DeepONet trunk basis (the smaller DeepONet), the error increases significantly, resulting in worse performance than the proposed method. 
This highlights the influence of the trunk network’s depth and parameter count. 
However, our method is able to reduce the number of parameters while simultaneously improving accuracy.

\textbf{Adaptivity.}
For all examples, we observed that the learned FEM basis exhibits adaptivity. 
Particularly, \textbf{(1)} for the Fokker–Plank (Section \ref{sec_fk_example}), aggregation–diffusion (Section \ref{sec_aggregation}), Keller–Segel (Section \ref{sec_keller}), and KdV (Section \ref{sec_kdv}) equation, where we intentionally introduced localized structures (bumps or peaks), the learned FEM basis functions adaptively concentrate around these regions in order to accurately capture the local features.
Meanwhile, \textbf{(2)} for the Allen–Cahn problem (Section~\ref{sec_allen_cahn}) and the Cahn–Hilliard equation (Section~\ref{sec_cahn_hilliard}), the distinctive features of the solution landscapes, such as shocks, rapid growth, and rapid decay, are more evenly distributed across the domain among different training and testing samples. Consequently, although the hat basis functions $h$ are initialized with relatively large (global) supports, they shrink to smaller supports after training, exhibiting localized behavior; 
at the same time, their centers remain approximately evenly distributed after learning, enabling the model to capture these features effectively throughout the domain.
These results show that our method effectively mimics the adaptive FEM and numerically justify its accuracy and efficiency.

\subsection{Allen-Cahn Equation}
\label{sec_allen_cahn}
In this section, we consider the Allen-Cahn equation,
\begin{align}
    u_t - {\epsilon^2} u_{xx}+ F'(u) = 0, x\in[0, 1], t\in[0, 10],
\end{align}
where $F(u) = (u^2-1)^2/4$, and $\epsilon = 0.01$. The target operator is the mapping from the initial condition to the solution at the terminal time. We generate the initial condition as $\lambda  \sin(2\pi x) + (1-\lambda)\sin(6 \pi (x -0.5 + \mu))$, where $\lambda$ and $\mu$ are free parameters. We will consider two sets of experiments with different free parameters settings.

\subsubsection{One degree of freedom}
\begin{figure}[H]
    \centering
    \includegraphics[scale = 0.5]{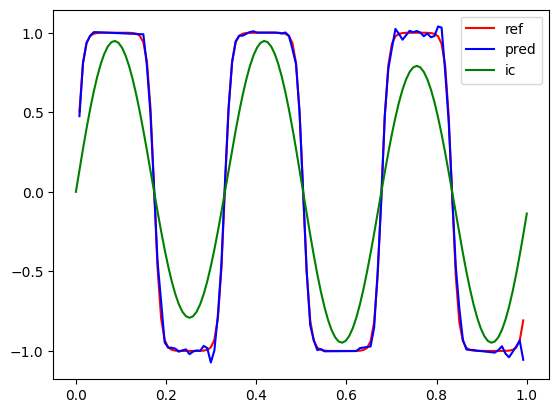}
    \includegraphics[scale = 0.5]{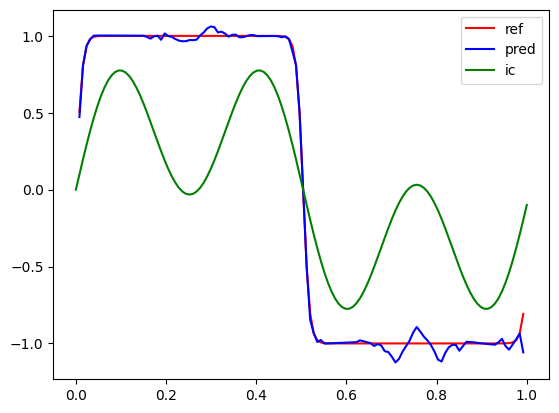}
    \caption{Demonstration of the predictions for the Allen Cahn example (one degree of freedom). The relative prediction error for left and right examples are $3.70\%$ and $4.00\%$ respectively. The average and variance of the relative error for 190 testing samples are $3.61\%\pm 1.3\%$.}
    \label{fig_ac_sol_pred}
\end{figure}

In this set of experiments, we fix $\mu = 0.5$ and sample $\lambda$ uniformly from $[0, 1]$.
We plot two typical solutions in Figure \ref{fig_ac_sol_pred}, and present the results in Table \ref{tab_alan_cahn}.

\begin{table}[H]
\centering
\begin{tabular}{|c|c|c|c|c|}
\hline
 & DeepONet-S & DeepONet-L& FERN & POD (same mesh) \\
\hline
Relative $L^2$ error         & $55.76\% \pm 25.77
\%$ & $3.39\% \pm 1.94\%$ & $3.61\% \pm 1.3\%$ & $0.15\% \pm 0.17\%$ \\
\hline
$\#$ Coefficient Parameters  & 19,200 & 19,200 & 19,200 & 19,200 \\
\hline
$\#$ Basis Parameters        & 4,200  & 54,700 & 80 & 0 \\
\hline
$\#$ Total Parameters        & 23,400 &73,900 & 19,280 & 19,200 \\
\hline
\end{tabular}
\caption{Comparison of different models for the Allen–Cahn example.
All models use a shared branch network architecture ($\hat{c}_k$) and employ 40 basis functions. For DeepONet, two configurations are considered: (1) DeepONet-S, which uses 2 fully connected layers to construct the basis functions, and (2) DeepONet-L, which uses 7 layers. 
The POD basis achieves the best accuracy among all methods. However, its use requires that both training and testing be performed on the same mesh, and that all training output functions be evaluated on this mesh. In contrast, the other models are evaluated on a denser mesh than used during training, demonstrating greater flexibility.
Additionally, the performance of the POD method drops significantly in the following examples, where more degrees of freedom are used in data generation, resulting in a higher-rank dataset.
}
\label{tab_alan_cahn}
\end{table}

\begin{figure}[H]
    \centering
    \includegraphics[trim={0 0 0 0.8cm},clip,scale = 0.5]{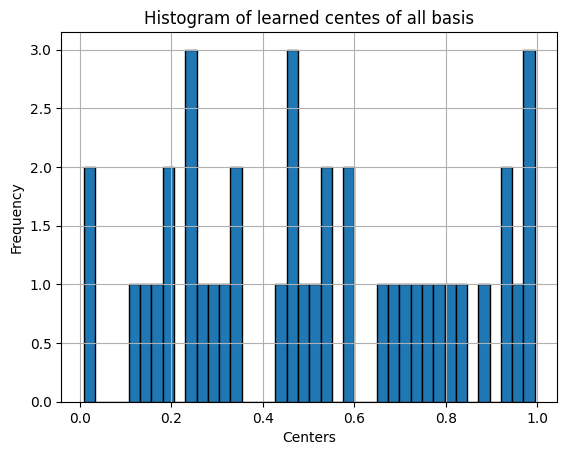}
    \includegraphics[trim={0 0 0 0.8cm},clip,scale = 0.5]{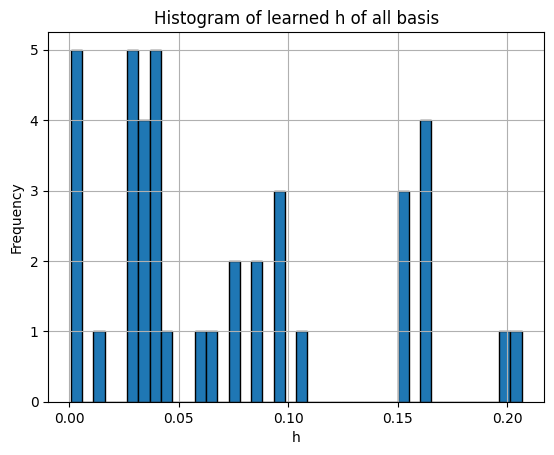}
    \caption{Allen-Cahn equation with one degree of freedom. Left: Histogram of all learned FEM basis centers $a_k$, which are initially uniformly distributed over $[0, 1]$. Right: Histogram of all learned support sizes $h_k$, all of which are initialized to $0.05$.
    Since the distinctive features of the solution landscapes, such as shocks, bumps, rapid decay, and rapid growth, are evenly distributed across the domain for different samples (see Figure \ref{fig_ac_sol_pred}), most learned hat functions shrink to smaller supports after training (right panel), exhibiting localized behavior. 
    Meanwhile, their centers remain evenly distributed to effectively capture these features across various samples (left panel). 
    This behavior is consistently observed for other numbers of basis functions tested; see Figure~\ref{fig_ac_2dofs_h_moving} for an illustration.
    }
    \label{fig_ac_hist_center_h}
\end{figure}

\textbf{Setting details.}
We generate a total of 167 input functions (initial conditions, ICs) along with their corresponding solutions for training purposes. For each output function, 100 evaluation points are uniformly sampled from the spatial domain. Each input function is discretized using 22 uniformly spaced sensors, corresponding to a mesh size of 22.
The model is trained using a cosine annealing learning rate schedule.
All neural network architecture are constructed using 40 basis functions (DeepONet global basis, FEM basis or POD basis). Each coefficient branch network follows a fully connected architecture of size $22 \times 20 \rightarrow 20 \times 1$, with the hyperbolic tangent (Tanh) function as the activation. The centers of the FEM basis functions are initialized uniformly over the interval $[0, 1]$, with a fixed support size $h = 0.05$.
A histogram of the learned FEM basis centers and $h$ after training is presented in Figure~\ref{fig_ac_hist_center_h}.

\textbf{Results analysis.}
As shown in Table~\ref{tab_alan_cahn}, the proposed FERN method achieves results comparable to the standard DeepONet while using significantly fewer trainable parameters. 
This demonstrates the effectiveness of the proposed approach: when a shallow network is designed to construct the FEM basis, a limited number of parameters is sufficient to achieve low relative prediction error. 

To further verify this observation, we also decrease the number of layers used in the trunk (basis) network. 
As shown in Table~\ref{tab_alan_cahn}, when the number of fully connected layers is reduced to two, the prediction error increases significantly. 
Although the POD method yields the smallest error, it relies on stronger assumptions when performing POD on the output functions, for instance, all outputs must be discretized on the same mesh. 
This restriction makes the POD-based approach a fixed-grid-to-fixed-grid mapping, which effectively corresponds to learning a function rather than an operator. Moreover, when the number of free parameters increases to two, the discretized output functions exhibit higher rank, causing the POD method to produce a larger error ($7.76\%$), which exceeds the $4.44\%$ error obtained using the proposed FEM basis. See the next section and Table~\ref{tab_alan_cahn_2dofs} for further details.

\subsubsection{Two degrees of freedom}
\begin{figure}[H]
    \centering
    \includegraphics[scale = 0.5]{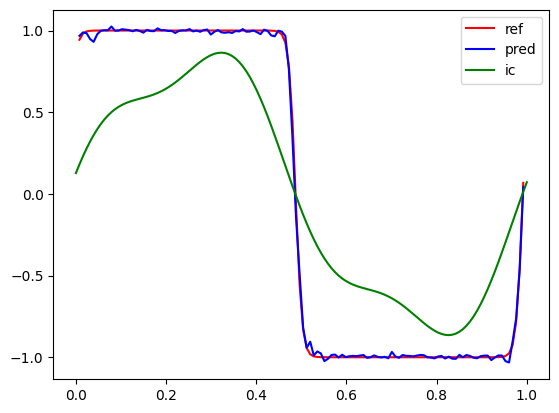}
    \includegraphics[scale = 0.5]{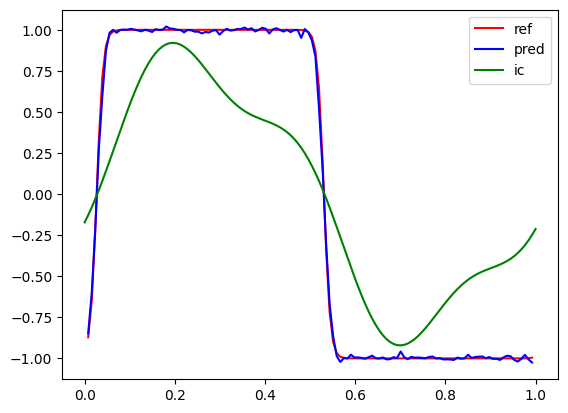}
    \includegraphics[scale = 0.5]{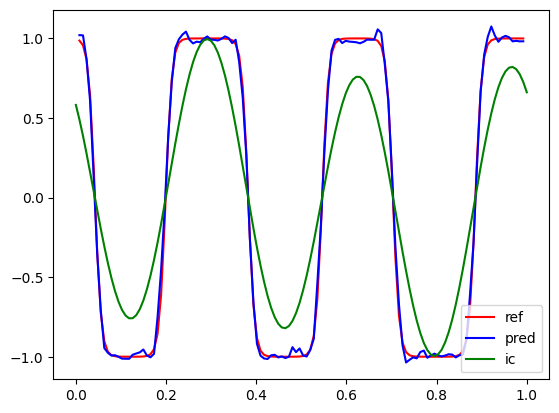}
    \includegraphics[scale = 0.5]{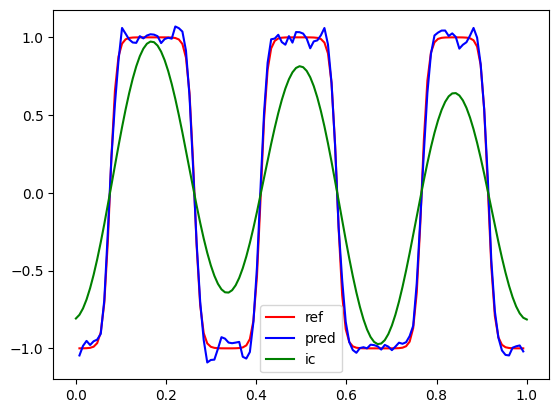}
    \caption{Predicted solutions for the Allen–Cahn equation with two degrees of freedom. The relative prediction errors for the four cases (displayed from left to right, top to bottom) are $2.21\%$, $2.23\%$, $3.90\%$, and $4.76\%$, respectively.}
    \label{fig_ac_2dofs_sol_pred}
\end{figure}

In this set of experiments, $\lambda$ and $\mu$ are both sampled from $[0, 1]$ uniformly; as a results, the solutions present more variability and make the problem more challenging. We plot four typical solutions in Figure \ref{fig_ac_2dofs_sol_pred}, and the results in Table \ref{tab_alan_cahn_2dofs}.

\begin{table}[H]
\centering
\begin{tabular}{|c|c|c|c|}
\hline
 & DeepONet& FERN & POD (same mesh) \\
\hline
Relative $L^2$ error         &  $8.65\% \pm 4.08\%$ & $4.44\% \pm 2.2\%$ & $7.76\% \pm 5.52\%$ \\
\hline
$\#$ Coefficient Parameters & 38,400 & 38,400 & 38,400 \\
\hline
$\#$ Basis Parameters        & 58,700 & 160 & 0 \\
\hline
$\#$ Total Parameters       &97,100 & 38,560 & 38,400 \\
\hline
\end{tabular}
\caption{Comparison of different models for the Allen–Cahn example with two degrees of freedom. 
All models employ the same branch networks architecture ($\hat{c}_k$) and approximate the output functions using 80 basis functions. 
We additionally test DeepONet with a two-layer trunk network; although this reduces the parameter count to $46,600$, the relative error stabilizes at $63.25\%$, which is substantially larger than that obtained with FERN.
Additionally, compared to the previous examples where a single degree of freedom was used to generate the training data, the accuracy of the POD method drops as the rank of the dataset increases.
}
\label{tab_alan_cahn_2dofs}
\end{table}

\begin{figure}[H]
    \centering
    \includegraphics[trim={0 0 0 0.8cm},clip,scale = 0.5]{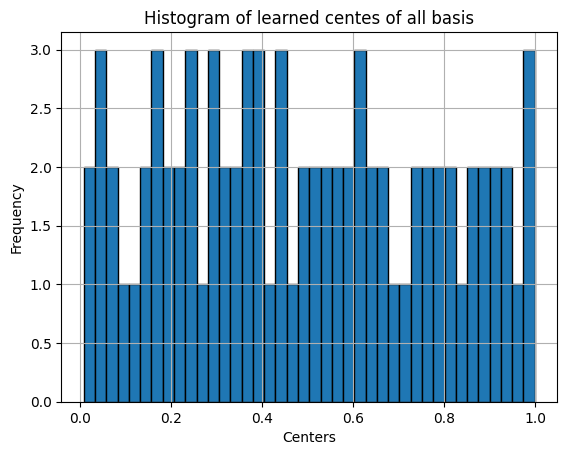}
    \includegraphics[trim={0 0 0 0.8cm},clip,scale = 0.5]{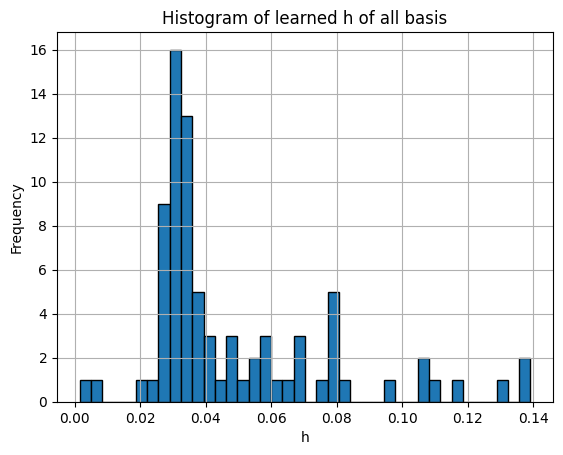}
    \caption{Allen-Cahn equation with two degrees of freedom. Left: Histogram of all learned FEM basis centers $a$, which are initially uniformly distributed over $[0, 1]$. Right: Histogram of all learned support sizes $h$, all of which are initialized to $0.05$.
    Since the distinctive features of the solution landscapes, such as shocks, bumps, rapid decay, and rapid growth, are evenly distributed across the domain for different samples (see Figure \ref{fig_ac_2dofs_sol_pred}), most learned functions shrink to smaller supports after training (right picture), exhibiting localized behavior.
    However, their centers remain evenly distributed to effectively capture these features across various samples (left picture).
   This behavior is consistently observed across different numbers of basis functions $N$ tested (see Figure~\ref{fig_ac_2dofs_h_moving}).
    }
    \label{fig_ac_2dofs_hist_center_h}
\end{figure}

\textbf{Setting details.}
We generate a total of 250 input functions (initial conditions, ICs) along with their corresponding solutions for training. Each output function is evaluated at 100 uniformly sampled points from the spatial domain, while each input function is discretized using 22 uniformly spaced sensors, corresponding to a mesh size of 22.
All models are trained for 2,000 epochs using a cosine annealing learning rate schedule. The neural network architectures across all models employ 80 basis functions, either the proposed FEM basis, the POD basis, or the standard DeepONet basis. Each coefficient branch network adopts a fully connected architecture with layer sizes $22 \times 20 \rightarrow 20 \times 1$.
For DeepONet, we report results using the ReLU activation function, which yielded a better accuracy compared to Tanh activation. 
For the FEM and POD models, the coefficient networks use the Tanh activation function for better accuracies. 
Notably, FERN exhibits robustness to the choice of activation: using ReLU, it achieves a relative prediction error of $4.45\% \pm 2.0\%$.
The centers of the FEM basis functions are initialized uniformly over the interval $[0, 1]$, with a fixed support size $h = 0.05$.
A histogram of the learned FEM basis centers and $h$ after training is presented in Figure~\ref{fig_ac_2dofs_hist_center_h}.

\textbf{Analysis of results.}
From Table~\ref{tab_alan_cahn_2dofs}, it can be observed that the proposed method achieves the highest prediction accuracy with the smallest variance. Although the POD non-trainable basis uses slightly fewer trainable parameters ($38{,}400$ vs.~$38{,}560$), it fails to handle output discretizations on different meshes and yields a larger error due to the increased number of degrees of freedom used to generate the input functions. 
These results further demonstrate the accuracy and efficiency of the proposed method.

\begin{figure}[H]
    \centering
    \includegraphics[scale = 0.5]{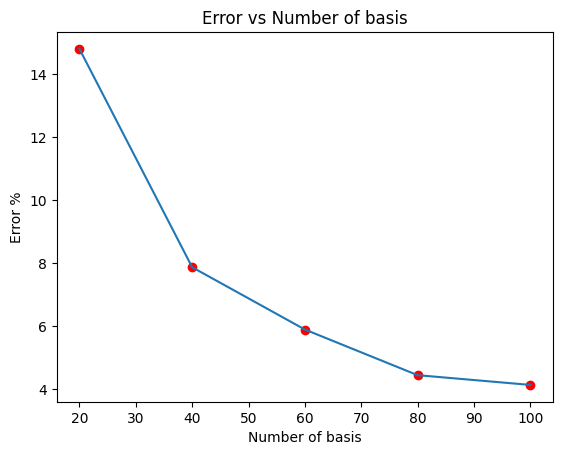}
    \includegraphics[scale = 0.5]{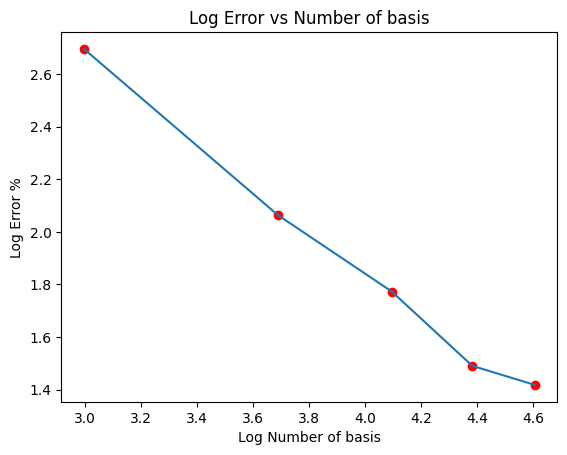}
    \caption{Allen-Cahn equation with two degrees of freedom and various number of basis. Left: Error decay with respect to the number of learned hat function basis. Right: decay in the log scale.}
    \label{fig_ac_error_decay}
\end{figure}

\begin{figure}[H]
    \centering
    \includegraphics[trim={0 0 0 0.8cm},clip,scale = 0.35]{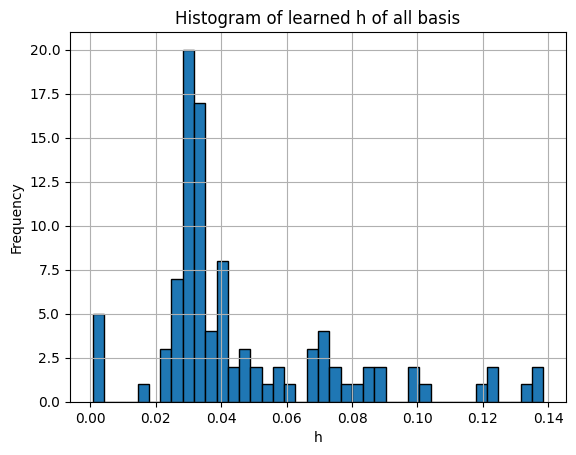}
    \includegraphics[trim={0 0 0 0.8cm},clip,scale = 0.35]{graphs/ac_2dofs_hist_h.png}
    \includegraphics[trim={0 0 0 0.8cm},clip,scale = 0.35]{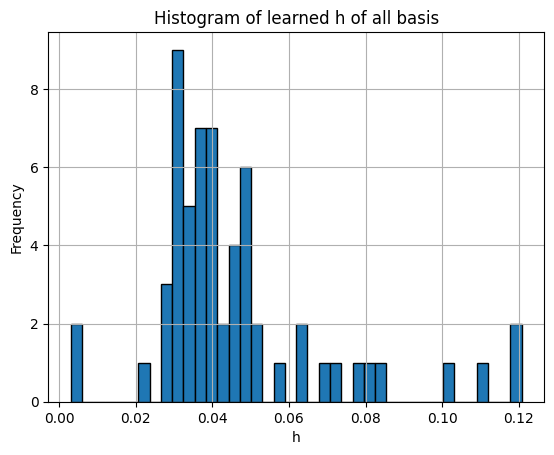}\\
    \includegraphics[trim={0 0 0 0.8cm},clip,scale = 0.35]{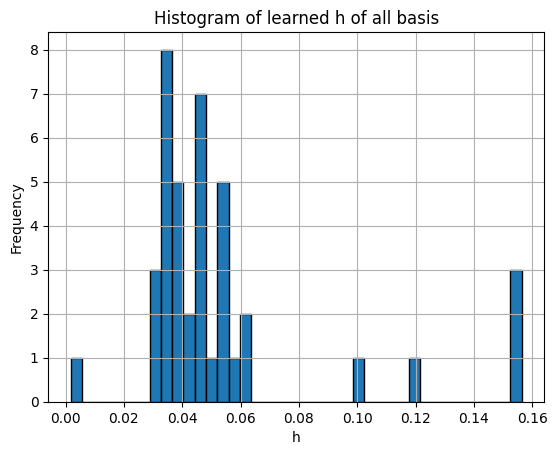}
    \includegraphics[trim={0 0 0 0.8cm},clip,scale = 0.35]{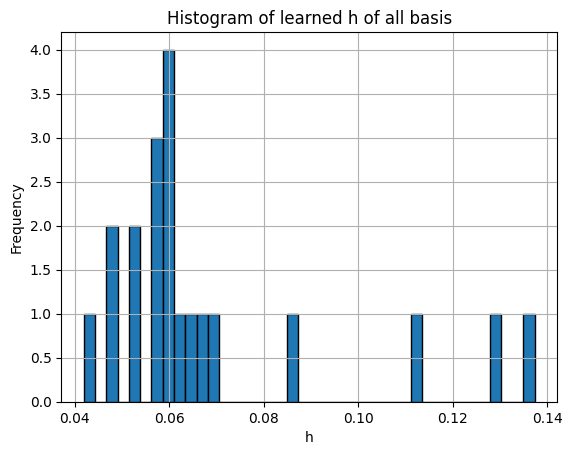}
    \caption{Allen-Cahn equation with two degrees of freedom and various number of basis. Distributions of the learned $h$ for models with different numbers of basis functions ($20$, $40$, $60$, $80$, and $100$). Each experiment is initialized with $h = 0.05$. As shown, the learned $h$ values shift across all cases, demonstrating the adaptivity of the proposed method.
    }
    \label{fig_ac_2dofs_h_moving}
\end{figure}

\textbf{A study on the number of basis.}
We investigate the relationship between prediction error and the number of basis networks (corresponding to the mesh size in the output space) for the proposed FEM-based methods. Specifically, we vary the number of learned FEM basis functions from 20 to 100 and train a separate model for each setting. 
All models are initialized with uniformly distributed centers in the domain [0, 1], and each basis function is assigned a fixed size of $h = 0.05$. The resulting error decay behavior is presented in Figure \ref{fig_ac_error_decay}.
Roughly, first order convergence is observed with respect to the number of basis.
In addition, we present the learned mesh size distributions for models with 20, 40, 60, 80, and 100 basis functions in Figure \ref{fig_ac_2dofs_h_moving}. From the figure, it can be observed that the learned mesh sizes gradually shift toward zero as the number of basis functions decreases. This indicates that, when fewer basis functions are used, the model tends to learn basis functions with larger support.

\subsection{Canh-Hilliard Equation}
\label{sec_cahn_hilliard}
In this section, we study the Canh-Hilliad equation,
\begin{align*}
    \frac{\partial u}{\partial t} = \frac{\partial}{\partial x} \left( M \frac{\partial}{\partial x} \left( -\epsilon^2 \frac{\partial^2 u}{\partial x^2} + F(u) \right) \right), x\in[0, 1], t\in[0, 10],
\end{align*}
where $F(u) = (u^2-1)^2/4$, and $M = 0.01$. The target operator is the mapping from the initial condition to the solution at the terminal time. We generate the initial condition by sampling free parameter $\lambda$ and $\mu$ for $\lambda  \sin(2\pi x) + (1-\lambda)\sin(6 \pi (x -0.5 + \mu))$. We will consider two sets of experiments with different free parameters settings.

\subsubsection{One degree of freedom}
\begin{figure}[H]
    \centering
    \includegraphics[scale = 0.5]{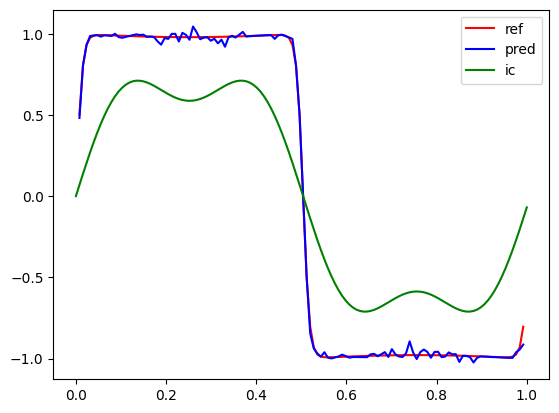}
    \includegraphics[scale = 0.5]{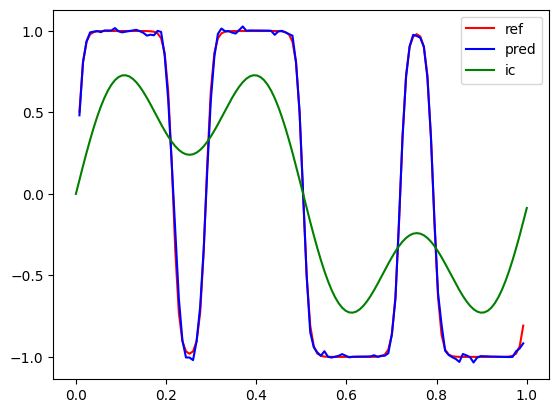}
    \caption{Demonstration of the predictions for the Cahn-Hilliard equation with one degree of freedom. The relative prediction error for left and right examples are $2.62\%$ and $2.19\%$ respectively. The average relative error for 190 testing samples is $2.62\%\pm 0.41\%$.}
    \label{fig_ch_1dof}
\end{figure}

In this set of experiments, we set the $\mu = 0.5$, and sample $\lambda$ uniformly from $[0, 1]$. We display two typical solutions in Figure \ref{fig_ch_1dof} and the numerical results in Table \ref{tab_ch_1dof}.

\begin{table}[H]
\centering
\begin{tabular}{|c|c|c|c|}
\hline
 & DeepONet& FERN & POD \\
\hline
Relative $L^2$ error        &  $2.91\% \pm 1.07\%$ & $2.62\% \pm 0.41\%$ & $0.38\% \pm 0.24\%$ \\
\hline
$\#$ Coefficient Parameters & 28,800  & 28,800 & 28,800 \\
\hline
$\#$ Basis Parameters        & 56,700  & 120 & 0 \\
\hline
$\#$ Total Parameters       & 85,500 & 28920 & 28,800 \\
\hline
\end{tabular}
\caption{Comparison of different models for the Cahn-Hilliard equation with one degree of freedom. All models share the same branch network architecture and have 60 basis functions (DeepONet global learnable basis, FEM learnable local basis, and POD basis). 
We also test a smaller DeepONet structure with a two-layer trunk network while keeping the branch network unchanged. Although this reduces the number of parameters to $35,000$, the relative error increases to $28.44\% \pm 11.07\%$.
}
\label{tab_ch_1dof}
\end{table}

\begin{figure}[H]
    \centering
    \includegraphics[trim={0 0 0 0.8cm},clip,scale = 0.5]{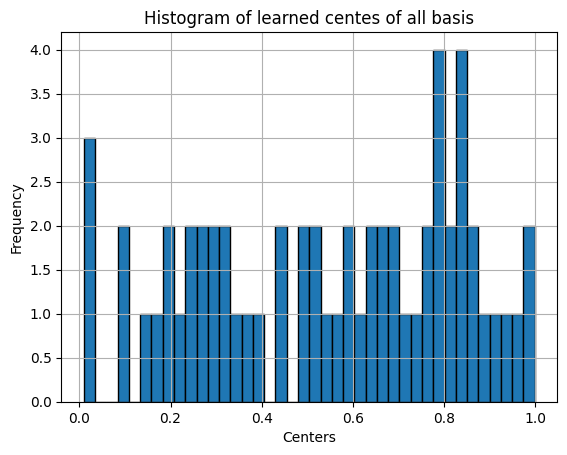}
    \includegraphics[trim={0 0 0 0.8cm},clip,scale = 0.5]{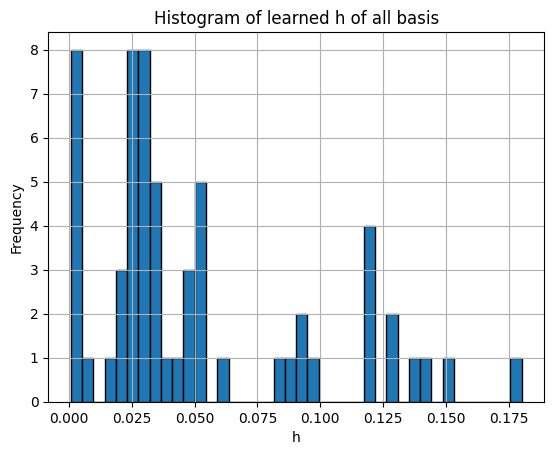}
    \caption{Cahn-Hilliard equation with one degree of freedom. Left: Histogram of all learned FEM basis centers, which are initially uniformly distributed over $[0, 1]$. Right: Histogram of all learned support sizes $h$, all of which are initialized to $0.05$.
    Since the distinctive features of the solution landscapes, such as shocks, bumps, rapid decay, and rapid growth, are evenly distributed across the domain for different samples (see Figure \ref{fig_ch_1dof}), most learned functions shrink to smaller supports ($0.025$) after training (left picture), exhibiting localized behavior.
    However, their centers remain almost evenly distributed to effectively capture these features across various samples (right picture).
    }
    \label{fig_ch_1dof_hist_center_h}
\end{figure}

\textbf{Setting details.}
We generate a total of 250 input functions (initial conditions, ICs) along with their corresponding solutions for training. Each output function is evaluated at 100 uniformly sampled points from the spatial domain, while each input function is discretized using 22 uniformly spaced sensors, corresponding to a mesh size of 22.
All models are trained for 2,000 epochs using a cosine annealing learning rate schedule. The neural network architectures across all models employ 60 basis functions—either the proposed FEM basis, the POD basis, or the standard DeepONet basis.
Each coefficient branch network adopts a fully connected architecture with layer sizes $22 \times 20 \rightarrow 20 \times 1$ and are activated by ReLU.
For the DeepONet and FEM models, the basis networks use the Tanh activation function.
The centers of the FEM basis functions are initialized uniformly over the interval $[0, 1]$, with a fixed support size $h = 0.05$.
A histogram of the learned FEM basis centers and $h$ after training is presented in Figure \ref{fig_ch_1dof_hist_center_h}.

\textbf{Analysis of the results.}
From Table~\ref{tab_ch_1dof}, the proposed method achieves a lower error ($2.62\% < 2.91\%$) compared to the standard DeepONet with a deep network-constructed basis. 
At the same time, the number of trainable parameters is significantly smaller, as we use only $2 \times 60$ parameters (two parameters per basis function) to construct a shallow network for building the FEM basis. 
Although the POD method outperforms the proposed method in terms of accuracy, it relies on non-trainable POD bases, requiring strong assumptions before application. For instance, all output functions must be discretized on the same mesh, which greatly limits the flexibility of the training data and restricts the POD-based operator learning to mappings from fixed grids to fixed grids, essentially resembling standard function learning rather than general operator learning.
Additionally, when we increase the number of degrees of freedom used to generate the training input functions, the output functions exhibit greater variability (i.e., higher rank). 
Consequently, the accuracy of the POD method drops significantly, becoming worse than that of the proposed method. See the results of the next experiment in Table~\ref{tab_ch_2dofs} for details.

Lastly, we observe the basis localization during network training. 
As shown in Figure \ref{fig_ch_1dof_hist_center_h}, all FEM hat functions are initialized with a uniform width of $h_k = 0.05$. However, after training, most functions shrink to narrower supports, exhibiting localized behavior. Since the special features of the output functions—such as sharp decays or mild shocks (see Figure \ref{fig_ch_1dof} for illustration)—are relatively evenly distributed across the domain, the centers $a_k$ of the learned functions do not shift significantly. 
In contrast, for examples where shocks are concentrated in specific regions of the domain, we observe clear localization of the functions around those regions.
See FK, KS, AD, KdV sections for examples.

\subsubsection{Two degrees of freedom}
\begin{figure}[H]
    \centering
    \includegraphics[scale = 0.5]{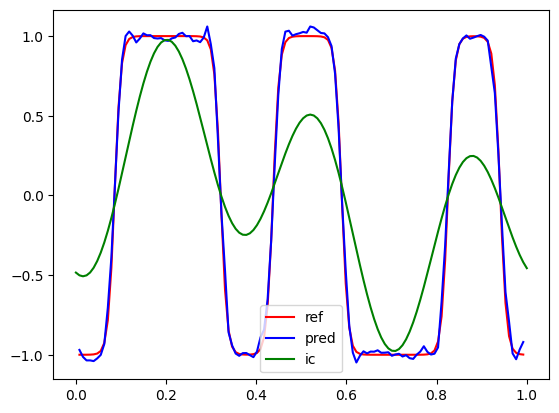}
    \includegraphics[scale = 0.5]{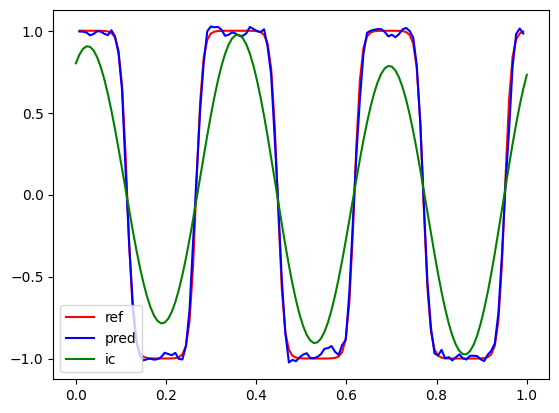}
    \includegraphics[scale = 0.5]{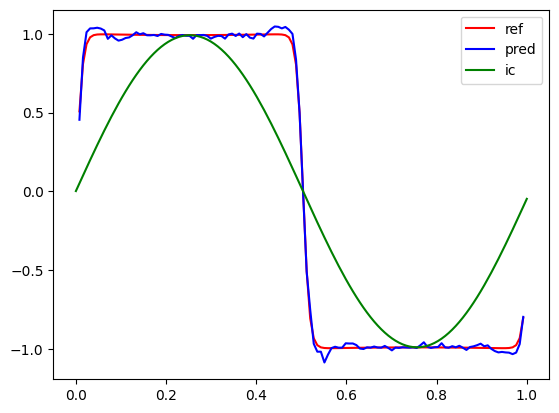}
    \includegraphics[scale = 0.5]{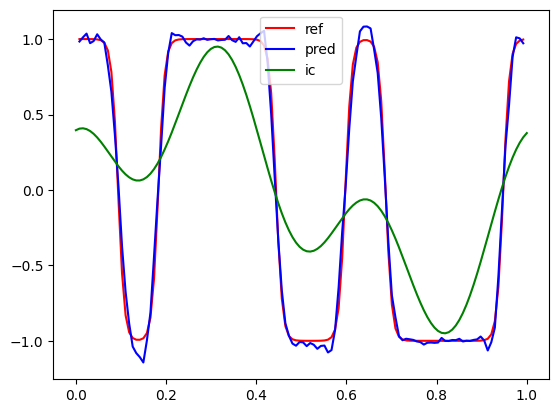}
    \caption{Predicted solutions for the Cahn-Hilliard equation with two degrees of freedom. The relative prediction errors for the four cases (displayed from left to right, top to bottom) are $4.09\%$, $4.29\%$, $2.59\%$, and $6.89\%$, respectively.}
    \label{fig_ch_2dofs_sols}
\end{figure}

In this set of experiments, we set $\lambda$ and $\mu$ both are sample uniformly from $[0, 1]$.
We plot four typical solutions in Figure \ref{fig_ch_2dofs_sols} and present the numerical results in Table \ref{tab_ch_2dofs}.
\begin{table}[H]
\centering
\begin{tabular}{|c|c|c|c|}
\hline
 & DeepONet& FERN & POD \\
\hline
Relative $L^2$ error         &  $14.68\% \pm 4.19\%$ & $5.59\% \pm 2.2\%$ & $7.22\% \pm 3.41\%$ \\
\hline
$\#$ Coefficient Parameters & 120,000 & 120,000 & 61,440 \\
\hline
$\#$ Basis Parameters        & 75,700  & 500 & 0 \\
\hline
$\#$ Total Parameters       & 195,700 & 120,500 & 61,440 \\
\hline
\end{tabular}
\caption{Comparison of different models for the Cahn-Hilliard equation with two degrees of freedom. All models share the coefficient networks; the DeepONet and FERN network have 250 learnable basis, while the POD network has 128 basis the largest number of basis due to the discretization size.
}
\label{tab_ch_2dofs}
\end{table}

\begin{figure}[H]
    \centering
    \includegraphics[trim={0 0 0 0.8cm},clip,scale = 0.5]{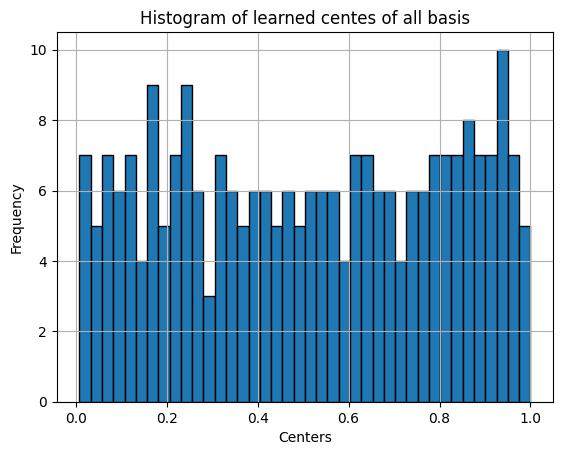}
    \includegraphics[trim={0 0 0 0.8cm},clip,scale = 0.5]{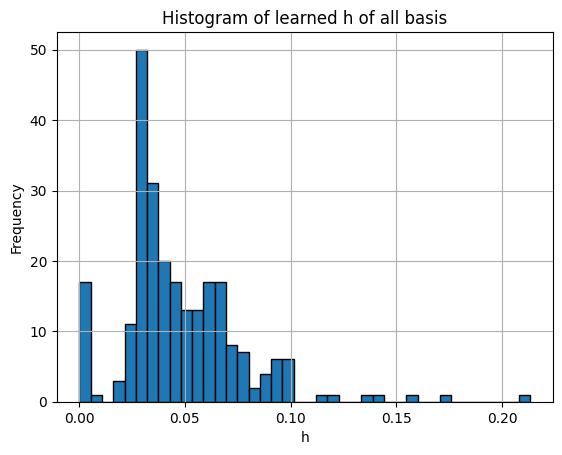}
    \caption{Cahn-Hilliard equation with two degrees of freedom. Left: Histogram of all learned FEM basis centers, which are initially uniformly distributed over $[0, 1]$. Right: Histogram of all learned support sizes $h$, all of which are initialized to $0.05$. 
    Most learned FEM hat functions shrink to smaller supports after training (left picture), exhibiting localized behavior.
    However, their centers remain evenly distributed to effectively capture these features across various samples (right picture).
    }
    \label{fig_ch_2dofs_hist_center_h}
\end{figure}

\textbf{Setting details.}
We generate a total of 1,000 input functions (initial conditions, ICs) along with their corresponding solutions for training. Each output function is evaluated at 100 uniformly sampled points from the spatial domain, while each input function is discretized using 22 uniformly spaced sensors, corresponding to a mesh size of 22.
All models are trained for 2,000 epochs using a cosine annealing learning rate schedule. The neural network architectures across all models employ 250 basis functions, either the proposed FEM basis, the POD basis, or the standard DeepONet basis.
Each coefficient branch network adopts a fully connected architecture with layer sizes $22 \times 20 \rightarrow 20 \times 1$ and are activated by ReLU.
For the DeepONet and FEM models, the basis trunk networks use the Tanh activation function.
The centers of the FEM basis functions are initialized uniformly over the interval $[0, 1]$, with a fixed support size $h = 0.05$.
A histogram of the learned FEM basis centers and $h$ after training is presented in Figure \ref{fig_ch_2dofs_hist_center_h}.

\textbf{Analysis of the results.} Similar to our observations for the Allen-Cahn equation, the proposed method achieves the ghighest prediction accuracy with the smallest variance.

\subsection{Fokker-Plank Equation}
\label{sec_fk_example}
\begin{figure}[H]
    \centering
    \includegraphics[scale = 0.5]{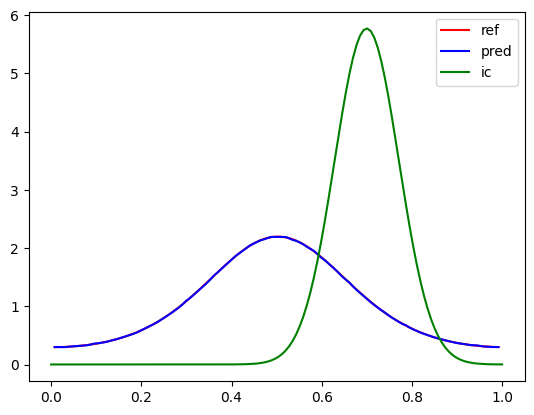}
    \includegraphics[scale = 0.5]{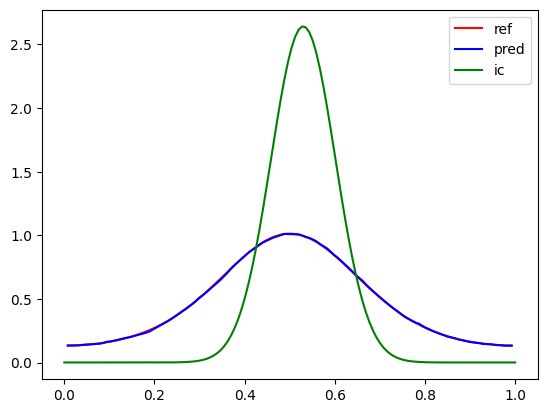}
    \caption{Fokker-Plank equation solution illustrations. The predictions are obtained using the proposed method with 30 basis functions, yielding errors of $0.28\%$ and $0.59\%$ respectively.}
    \label{fig_fk_sol}
\end{figure}

In this section, we consider the Fokker-Planck (FK) equation,
\begin{align}
    u_t - (u(\log(u) + \cos(2 \pi x))_x)_x = 0, x\in[0, 1], t\in[0, 0.1].
\end{align}
The target operator is the mapping from the initial condition to the solution at the terminal time. 
We generate the initial condition by sampling free parameters $c_0$ and $c_1$. Specifically, the initial condition has the form 
$
u_0(x)= c_1 \exp(-100(x-c_0)^2) + 10^{-3},
$
where $c_0$ is uniformly drawn from $[0.3, 0.7]$ and $c_1$ is from $[1, 10]$.
We plot two typical solutions in Figure \ref{fig_fk_sol}.
Notably, the terminal solutions exhibit a bell-shaped profile with a peak centered around 0.5 across different initial conditions. This consistent behavior provides an opportunity to evaluate the adaptivity of the learning process, that is, the model allocates more basis functions to the central region of the domain, where the solution is most concentrated, rather than to the flat tails.

We present the numerical results in Table~\ref{tab_fk_uniform} and Table~\ref{tab_fk_non_uniform}. Two sampling scenarios are considered:  
(1) In the first case, the terminal solutions corresponding to each initial condition are uniformly sampled over the entire domain with a mesh size of 64. This setting permits the use of the POD method to construct basis functions for prediction.  
(2) In the second case, one-third of the terminal solutions are sampled uniformly only over the interval \([0, 0.5]\), another third over \([0.5, 1]\), and the remaining third uniformly over the full interval \([0, 1]\). 
Under this non-uniform sampling scheme (different meshes for different output functions), POD cannot be applied directly, necessitating the use of learnable basis functions.
In both cases, our method achieves comparable or superior prediction accuracy while using significantly fewer trainable parameters. 
For example, with 30 basis functions, FERN requires only 14,460 trainable parameters, whereas the standard DeepONet uses 68,100.
\begin{table}[H]
\centering
\begin{tabular}{|p{2.2cm}|p{2.5cm}|p{2.5cm}|c|c|c|}
\hline
 & DeepONet (30) & 2-layer DeepONet (30) & FERN (30) & DeepONet (10) & FERN (10)  \\
\hline
Relative $L^2$ error         &  $1.16\% \pm 1.36\%$ &  $4.23\% \pm 3.49\%$ & {$1.10\% \pm 0.58\%$} & $1.28\% \pm 1.59\%$ & $1.15\% \pm 0.63\%$ \\
\hline
$\#$ Coefficient Parameters & 14,400 & 14,400 & 14,400 & 4,800 & 4,800   \\
\hline
$\#$ Basis Parameters        & 53,700 & 3,200  & 60 & 51,700 & 20  \\
\hline
$\#$ Total Parameters       & 68,100 & 17,600 & 14,460 & 56,500 & 4,820  \\
\hline 
\end{tabular}
\caption{Comparison of different models for the Fokker-Plank equation (the numbers $30$ and $10$ denote the number of basis $N$ used) on uniform grids. 
The terminal solutions are uniformly sampled over the interval $[0, 1]$ for all samples. Notably, all DeepONet models achieve accurate predictions when using a deep basis network; however, their performance deteriorates when the basis network is reduced to only two layers, see the third column. In contrast, FERN maintain high accuracy while requiring substantially fewer trainable parameters.
}
\label{tab_fk_uniform}
\end{table}

\begin{table}[H]
\centering
\begin{tabular}{|p{2.5cm}|c|p{2cm}|c|c|c|}
\hline
 & DeepONet (30) & 2-layer DeepONet (30) & FERN (30) & DeepONet (10) & FERN (10)  \\
\hline
Relative $L^2$ error         &  $1.08\% \pm 0.94\%$ &  $6.93\% \pm 4.34\%$ & $1.48\% \pm 0.9\%$ & $1.26\% \pm 1.27\%$ & $1.35\% \pm 0.83\%$ \\
\hline
$\#$ Coefficient Parameters & 14,400 & 14,400 & 14,400 & 4,800 & 4,800   \\
\hline
$\#$ Basis Parameters        & 53,700 & 3,200  & 60 & 51,700 & 20  \\
\hline
$\#$ Total Parameters        & 68,100 &  17,600 & 14,460 & 56,500 & 4,820   \\
\hline 
\end{tabular}
\caption{Comparison of different models for the Fokker-Plank equation (the numbers $30$ and $10$ denote the number of basis $N$ used) on non-uniform grids. In this setting, the terminal solutions are {not} uniformly sampled from the full domain \([0, 1]\) across all samples. Specifically, one-third of the samples are drawn from the interval \([0, 0.5]\), another one-third from \([0.5, 1]\), and the remaining one-third from the full interval \([0, 1]\). This non-uniform sampling pattern prevents the direct application of POD-based basis functions.
Similar to the uniform sampling results shown in Table~\ref{tab_fk_uniform}, the performance of the standard DeepONet deteriorates when the number of layers in the basis network is reduced to two, whereas FERN maintains accurate predictions while using significantly fewer trainable parameters.
}
\label{tab_fk_non_uniform}
\end{table}

\textbf{Detailed settings.}
We generate a total of 42 input functions (initial conditions, ICs) along with their corresponding solutions for training. 
Each input function is discretized using 22 uniformly spaced sensors, corresponding to a mesh of size 22. The corresponding output functions are evaluated at 64 spatial points, yielding a total of \(42 \times 64\) training samples.
Two different samples are tested, refer to Table \ref{tab_fk_uniform} and Table \ref{tab_fk_non_uniform} for details.
All FEM-based models are trained for 2,000 epochs using a cosine annealing learning rate scheduler, while DeepONet models are trained with more epochs to ensure convergence to a stable error level. 
For the standard DeepONet, the basis network follows a fully connected feedforward architecture with the structure
$1 \times 100 \rightarrow 100 \times 100 \rightarrow 100 \times 100 \rightarrow 100 \times 100 \rightarrow 100 \times 100 \rightarrow 100 \times 100 \rightarrow 100 \times K,$
where each hidden layer is followed by a \(\tanh\) activation, and the final layer is linear. 
All layers include bias terms except for the output layer. We also evaluate a variant of DeepONet with only two layers in the basis network, where the architecture is $1 \times 100 \rightarrow 100 \times 100$, and the activation function used is ReLU.

\textbf{Accuracy.}
On uniform grids, the performance of FERN is similar to that of DeepONet, but with a smaler variance, see Table \ref{tab_fk_uniform}. From Table \ref{tab_fk_non_uniform}, on nonuniform grids, when 30 basis functions are used, FERN achieves a prediction error of $1.48\%$, while the standard DeepONet with a deep trunk basis attains an error of $1.08\%$.
Although the proposed method does not yield a smaller average error in this example, it exhibits a smaller prediction variance ($0.9\%$ vs $0.94\%$) and requires significantly fewer trainable parameters ($14{,}460$ vs.~$68{,}100$), as it employs a two-layer shallow structure to construct the adaptive FEM basis.

To further demonstrate that not all shallow-structured basis networks are effective, we also reduce the number of fully connected layers in the standard DeepONet to two, which decreases the number of trainable parameters to $17{,}600$. 
However, despite having more parameters than the proposed method ($17{,}600$ vs.~$14{,}460$), the resulting error is considerably larger ($6.93\%$ vs.~$1.48\%$). This highlights the effectiveness of the learned adaptive FEM basis.

\begin{figure}[H]
    \centering
    \includegraphics[trim={0 0 0 0.8cm},clip,scale = 0.5]{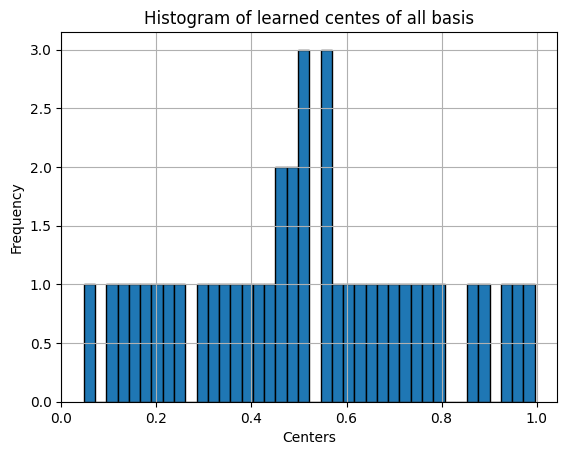}
    \includegraphics[scale = 0.5]{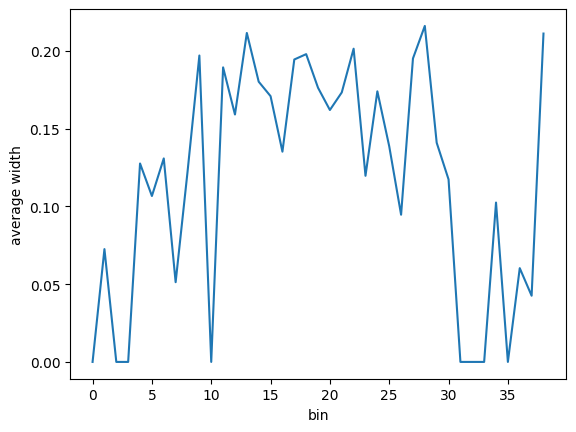}
    \caption{Fokker-Plank equation. The left panel shows the distribution of the learned basis function centers (40) across the domain. The right panel presents the average support width $h$ of the learned hat functions within each subinterval (bin), where the average is computed over all basis functions whose centers fall within the respective bin.
    We observe that, after training, more basis functions are concentrated near $x = 0.5$, where a bump appears in the solution landscape (see Figure~\ref{fig_fk_sol}), demonstrating the adaptivity of the learning.
    }
    \label{fig_fk_distribution_center}
\end{figure}

\textbf{Adaptivity.}
As shown in Figure \ref{fig_fk_sol}, all terminal solutions,represented in terms of the FEM basis, exhibit a bell-shaped profile with their peaks centered near the middle of the domain. Consequently, we expect the learned basis functions to reflect this spatial concentration pattern. This behavior is indeed observed in our results.
To further validate this hypothesis, we conduct experiments using a larger number of basis functions. In the left panel of Figure \ref{fig_fk_distribution_center}, we display the distribution of the learned FEM basis centers. 
The results indicate that a greater number of basis functions are concentrated near the center of the domain.
At the same time, the right panel shows the average support width $h$ of the learned hat functions, revealing that basis functions centered near the middle region of the domain tend to have wider supports compared to those located elsewhere.
Since hat functions are composed of two linear segments with slopes $1$ and $-1$, a larger support (width $h$) corresponds to a basis function of greater magnitude. This property may enhance the network's ability to represent the bell-shaped terminal solutions effectively.

\begin{figure}[H]
    \centering
    \includegraphics[trim={0 0 0 0.8cm},clip,scale = 0.35]{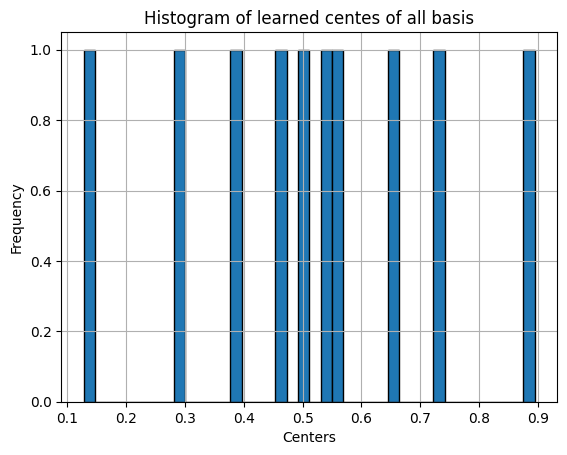}
    \includegraphics[trim={0 0 0 0.8cm},clip,scale = 0.35]{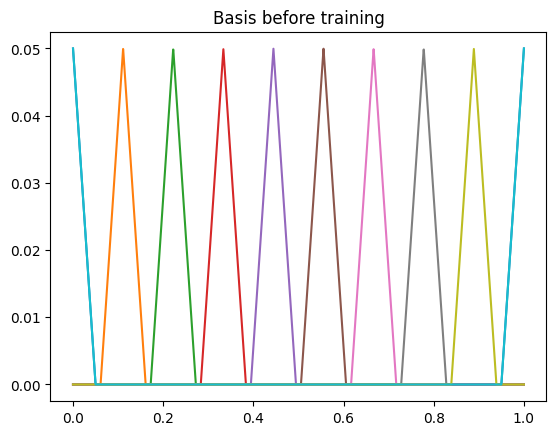}
    \includegraphics[trim={0 0 0 0.8cm},clip,scale = 0.35]{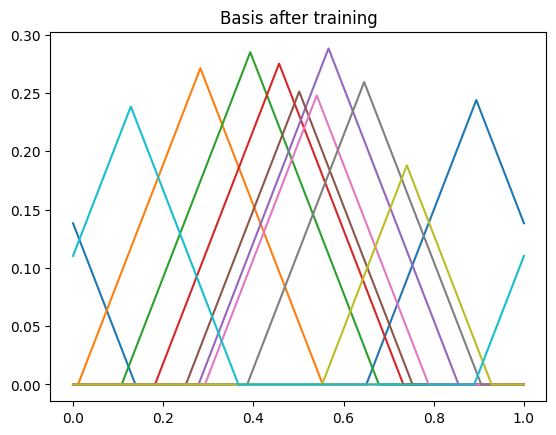}
    \caption{Fokker-Plank equation.
Demonstration of the basis functions (10 in total). Left: distribution of basis centers after training. Middle: all basis functions before training. Right: all basis functions after training.
We observe that, after training, more basis functions are concentrated near $x = 0.5$, where a bump appears in the solution landscape (see Figure~\ref{fig_fk_sol}), demonstrating the adaptivity of the learned FEM functions.
}
    \label{fig_fk_bs}
\end{figure}

We also present the basis functions before and after training for the case with 10 basis functions in Figure \ref{fig_fk_bs}. The initial centers are uniformly distributed across the domain with a fixed width $h = 0.05$.

\subsection{Aggregation-diffusion equation}
\label{sec_aggregation}
\begin{figure}[H]
    \centering
    \includegraphics[scale = 0.5]{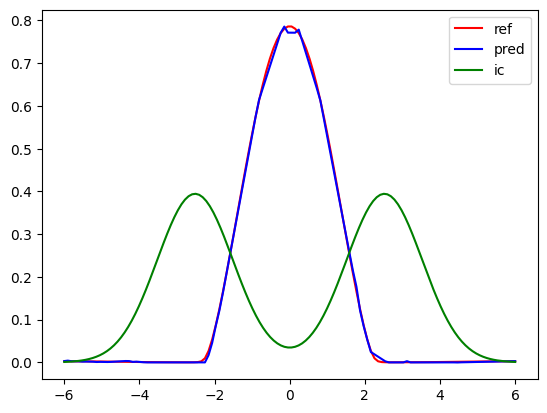}
    \includegraphics[scale = 0.5]{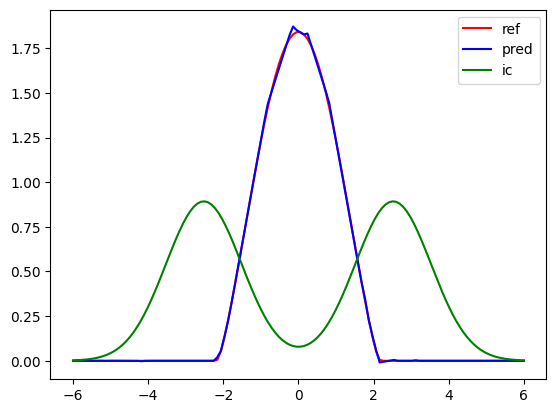}
    \caption{Demonstration of two solutions of the aggregation-diffusion equation.  Notably, the terminal solutions corresponding to different initial conditions remain close to zero outside the central region of the domain, indicating strong spatial localization.}
    \label{fig_aggre_diff_sol_demon}
\end{figure}

In this section, we study the aggregation-diffusion equation,
\begin{align*}
    u_t - (u(D u^{m-1} + W(x)*u)_x)_x = 0, \ \ x\in[-6, 6], \ \ t\in[0, 200],
\end{align*} 
where $D = 0.4, m = 2$, $W(x) =  -\frac{1}{\sqrt{2 \pi}\sigma} \exp\left(-\frac{x^2}{2\sigma^2}\right)$ with $\sigma = 1$, and $*$ is the convolution.
The target operator is the mapping from the initial
condition to the solution at the terminal time. We generate the initial condition by sampling the free parameter $c_0$ uniformly in $[1, 5]$ from $$u_0(x) = \frac{c_0}{2\sqrt{2\pi}}(\exp(-(x-x_0)^2/2)+\exp(-(x+x_0)^2/2)).$$
We present two representative solutions of the system in Figure \ref{fig_aggre_diff_sol_demon}. 
Notably, despite differing initial conditions, the terminal states consistently exhibit a bell-shaped profile centered within the domain, with the solution remaining nearly flat near the boundaries. This example serves as a test of the adaptability of the learned basis functions.
The coefficient networks (used in both the standard DeepONet and the FEM-based network) adopt a fully connected architecture with layer dimensions \(22 \times 20 \rightarrow 20 \times 1\), and all layers are activated by the hyperbolic tangent (\(\tanh\)) function.

We present the numerical results in Table \ref{tab_agg_diff_2dofs}.
The proposed methods can achieve similar accuracy while using fewer trainable parameters, given the same training and network settings. We give the numerical training details below.
\begin{table}[H]
\centering
\begin{tabular}{|c|c|c|c|c|}
\hline
 & DeepONet (40) & FERN (40) & DeepONet (20) & FERN (20)  \\
\hline
Relative $L^2$ error         &  $1.44\% \pm 2.69\%$ & $1.08\% \pm 0.54\%$ & $1.41\% \pm 2.74\%$ & $1.51\% \pm 2.39\%$ \\
\hline
$\#$ Coefficient Parameters & 19,200 & 19,200 & 9,600& 9,600   \\
\hline
$\#$ Basis Parameters        & 54,700  & 80 & 52,700 & 40  \\
\hline
$\#$ Total Parameters       & 73,900 & 19,280 & 62,300 & 9,640  \\
\hline 
\end{tabular}
\caption{Comparison of different models for the aggregation-diffusion equation (the numbers $40$ and $20$ denote the number of basis $N$ used). 
We also test a smaller DeepONet with a two-layer trunk network while keeping the branch network the same as in all other models. 
Although this reduces the total number of parameters to $23,400$ (still more than the proposed method with $19,280$ parameters), the relative error increases significantly to $61.63\% \pm 2.74\%$.
}
\label{tab_agg_diff_2dofs}
\end{table}
\textbf{Setting details.} 
We generate a total of 42 input functions (initial conditions, ICs) along with their corresponding solutions for training. Each input function is discretized using 22 uniformly spaced sensors, corresponding to a mesh of size 22. The corresponding output functions are evaluated at 64 uniformly sampled spatial points, yielding a total of \(42 \times 64\) training samples.
All models are trained for 2,000 epochs using a cosine annealing learning rate schedule.
We test two different settings for the number of basis functions, as summarized in Table~\ref{tab_agg_diff_2dofs}. In both settings, the coefficient network architecture (used in both the standard DeepONet and FERN) adopts a fully connected architecture with layer dimensions \(22 \times 20 \rightarrow 20 \times 1\), and all layers are activated by the hyperbolic tangent (\(\tanh\)) function.

The centers of the FEM basis functions are initialized uniformly over the input interval, with a fixed support width of \(h = 0.05\). For the standard DeepONet, the basis network follows a fully connected feedforward architecture with the structure
$1 \times 100 \rightarrow 100 \times 100 \rightarrow 100 \times 100 \rightarrow 100 \times 100 \rightarrow 100 \times 100 \rightarrow 100 \times 100 \rightarrow 100 \times K,$
where each hidden layer is followed by a \(\tanh\) activation, and the final layer is linear. All layers include bias terms except the output layer.

\textbf{Accuracy.}
As shown in Table~\ref{tab_agg_diff_2dofs}, the proposed FEM-basis method requires significantly fewer trainable parameters across both basis settings. 
When 20 basis functions are used, FERN attains a slightly higher average error than the standard DeepONet ($1.51\%$ vs.~$1.41\%$). 
However, when the number of basis functions increases to 40, FERN achieves a smaller error ($1.08\%$ vs.~$1.44\%$), a trend that is consistently observed in most of our experiments. 
Furthermore, even in cases where the average prediction error is higher, FERN consistently exhibits smaller variance in prediction errors, demonstrating its robustness and suggesting superior generalization performance in extrapolation scenarios.

\begin{figure}[H]
    \centering
    \includegraphics[scale = 0.5]{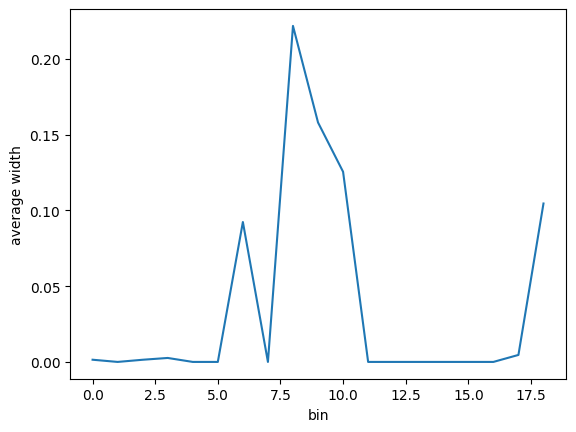}
    \includegraphics[scale = 0.5]{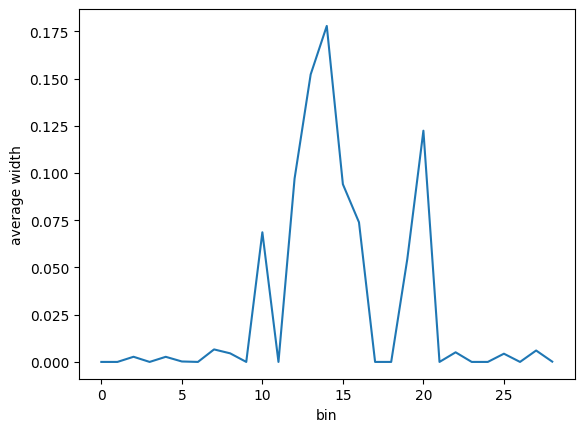}
    \includegraphics[scale = 0.5]{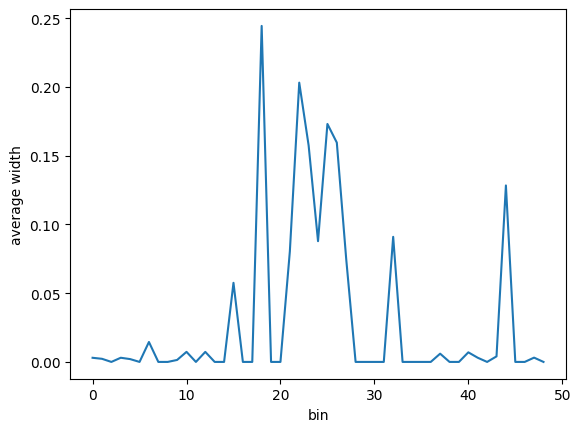}
    \includegraphics[scale = 0.5]{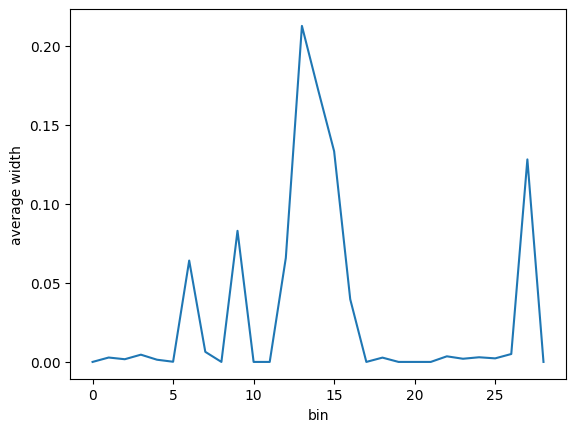}
    \caption{Aggregation-diffusion equation. Distribution of learned FEM basis support across the domain for different numbers of basis functions, ranging from 10 to 40 (ordered from left to right, top to bottom). The PDE domain is partitioned into uniform bins (subintervals), and for each bin, we compute the average support size \(h\) of the learned FEM basis functions whose centers $a_k$ fall within that bin. 
    The x-axis represents the bins from the left to the right boundary, while the y-axis indicates the corresponding average support.  
   Note that larger support values result in wider hat functions, since the slopes of the underlying piecewise linear components are fixed at 1 and -1.
 }
    \label{fig_aggre_diff_h_distribution}
\end{figure}

\textbf{Adaptivity.}
As shown in Figure \ref{fig_aggre_diff_sol_demon}, the terminal solutions exhibit a bell-shaped profile and decay rapidly outside the center of the domain. Since each hat basis function consists of two linear segments with slopes 1 and -1, a larger support width $h$ corresponds to a basis function with greater spatial extent. Therefore, we expect the learned hat functions to have larger support when their centers are located near the center of the domain, and smaller support as their centers move toward the boundaries. 
This behavior is confirmed in our numerical experiments. To validate this hypothesis, we vary the number of basis functions from 10 to 40 and plot the distribution of the learned support widths $h$ across uniformly divided subintervals of the domain. See Figure \ref{fig_aggre_diff_h_distribution} for the results.

\subsection{Keller-Segel Equation}
\label{sec_keller}
\begin{figure}[H]
    \centering
    \includegraphics[scale = 0.5]{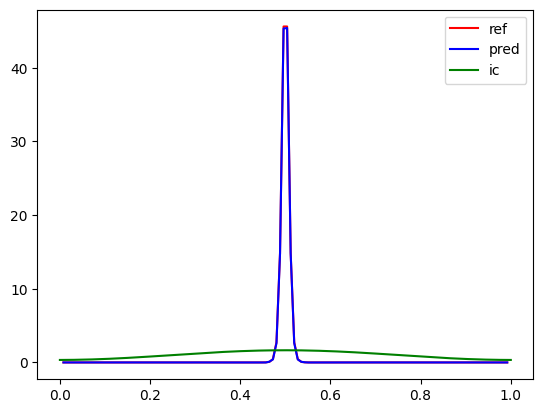}
    \includegraphics[scale = 0.5]{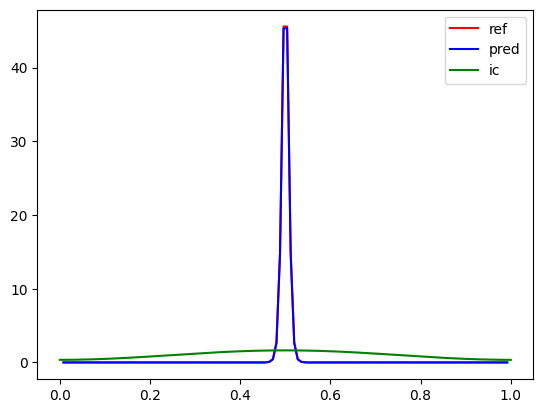}
    \caption{Demonstration of two solutions of the Keller-Segel equation.  Notably, the terminal solutions corresponding to different initial conditions exhibit a bump close to $x = 0.5$.}
    \label{fig_aggre_ks_sol_demon}
\end{figure}

In this section, we study the Keller-Segel (KS) equation,
\begin{align*}
\begin{cases}
    u_t - (u(D \log(u) - \chi v )_x)_x = 0,\\ 
    -v_{xx} + v = u, 
\end{cases}x\in[0, 1], t\in [0, 1.0]
\end{align*}
where $D = 0.01, \chi = 5.0$.
The target operator is the mapping from the initial
condition to the solution at the terminal time. 
We generate the initial condition by sampling the free parameter $c_0$ uniformly in $[0.2, 0.8]$ from $$u_0(x) = 1+c_0\sin(2\pi(x-0.25)).$$ 
We present two typical solutions of the system in Figure \ref{fig_aggre_ks_sol_demon}. 
Notably, despite differing initial conditions, the terminal states consistently exhibit a bump profile centered within the domain, with the solution remaining nearly flat near the boundaries. This example serves as a test of the adaptability of the learned basis functions to capture a sharp gradient in the solution landscape.

We present the numerical results in Table \ref{tab_ks_results}.
The proposed methods can achieve similar accuracy while using fewer trainable parameters, given the same training and network settings. We give the numerical training details below.
\begin{table}[H]
\centering
\begin{tabular}{|c|c|c|c|c|}
\hline
 & DeepONet(40) & 2-layer DeepONet & FERN (40) & FERN (60)  \\
\hline
Relative $L^2$ error         &  $0.73\% \pm 0.05\%$ & $5.97\% \pm 0.23\%$ & $1.08\% \pm 0.67\%$ & $0.78\% \pm 0.29\%$ \\
\hline
$\#$ Coefficient Parameters & 19,200 & 19,200 & 19,200 & 28,800   \\
\hline
$\#$ Basis Parameters        & 54,700  & 4,200 & 80 & 120  \\
\hline
$\#$ Total Parameters       & 73,900 & 23,400 & 19,280 & 28,920  \\
\hline 
\end{tabular}
\caption{Comparison of different models for the Keller-Segel equation (the numbers $40$ and $60$ denote the number of basis $N$ used). 
All proposed FEM-based models employ only a single layer to construct the FEM basis. Notably, when the number of basis layers in the standard DeepONet is reduced to two, its performance degrades. 
In contrast, FERN maintains high accuracy with a substantially smaller number of trainable parameters.
}
\label{tab_ks_results}
\end{table}

\textbf{Setting details.} 
We generate a total of 42 input functions (ICs) and their corresponding solutions for training, where each input function is discretized using 22 uniformly spaced sensors (i.e., a mesh of size 22). 
Different output functions are evaluated on different non-uniform meshes to introduce challenges for using POD bases: specifically, one-third of the output functions are sampled only on a fine mesh (size 49) in the first third of the domain $[0, 1]$, another third only in the middle third, and the final third only in the last third of the domain.
All FEM-based models are trained for 2,000 epochs using a cosine annealing learning rate schedule, while DeepONet models are trained with significantly more epochs to reach stable performance. 
We explore different configurations for the number of basis functions and the depth of the basis network in DeepONet, as summarized in Table \ref{tab_ks_results}. 

In all settings, the coefficient network architecture (used in both DeepONet and FERN) follows a fully connected architecture with dimensions $22 \times 20 \rightarrow 20 \times 1$, and use the hyperbolic tangent ($\tanh$) activation for all layers. 
The centers of the FEM basis functions are uniformly initialized across the domain with a fixed support width of $h = 0.05$. 
For the standard DeepONet, the basis network adopts a deep fully connected architecture with the structure $1 \times 100 \rightarrow 100 \times 100 \rightarrow 100 \times 100 \rightarrow 100 \times 100 \rightarrow 100 \times 100 \rightarrow 100 \times 100 \rightarrow 100 \times N$, where $N$ is the number of basis, with $\tanh$ activations following each hidden layer and a linear output layer; all layers include bias terms except the final one. For the two-layer DeepONet variant, the basis network consists of two layers with dimensions $1 \times 100 \rightarrow 100 \times 1$.

\textbf{Accuracy.}
Firstly, due to the use of a non-uniform mesh, as detailed in the previous section, the POD method cannot be applied. 
As shown in Table~\ref{tab_ks_results}, consistent with the results observed in other examples, the proposed method achieves comparable or even superior average relative prediction error while using substantially fewer trainable parameters.
In particular, when a standard DeepONet is constructed with a two-layer fully connected trunk network, the average error is $5.97\%$, whereas the proposed method with the learnable FEM basis achieves a significantly lower error of $1.08\%$. By increasing the depth of the basis network of DeepONet or the number of basis for the propose method, both accuracies can be improved. While the proposed method uses fewer number of parameters.

\begin{figure}[H]
    \centering
    \includegraphics[trim={0 0 0 0.8cm},clip,scale = 0.5]{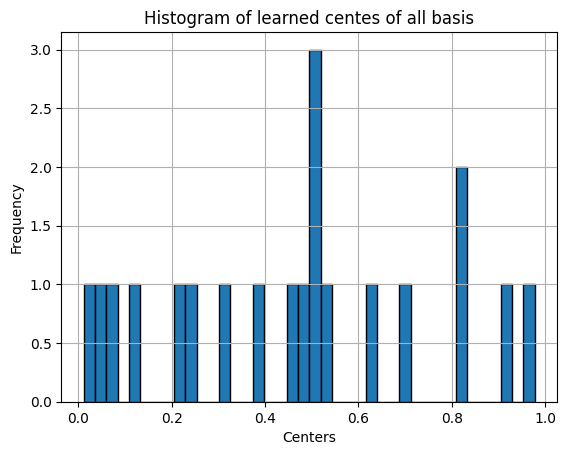}
    \includegraphics[trim={0 0 0 0.8cm},clip,scale = 0.5]{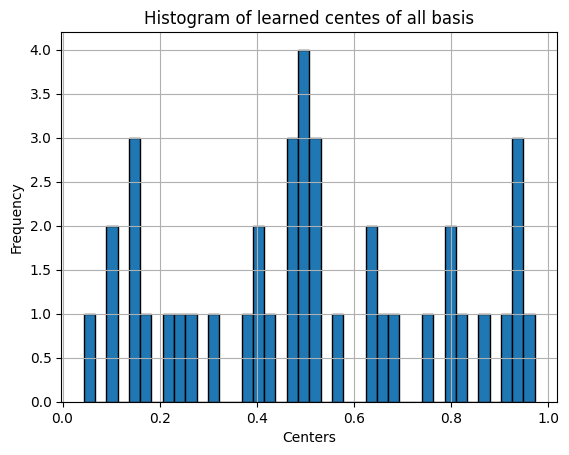}\\\vspace{0.5cm}
    \includegraphics[trim={0 0 0 0.8cm},clip,scale = 0.5]{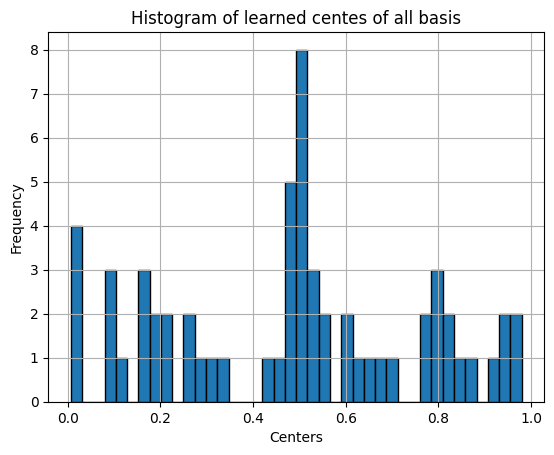}
    \includegraphics[trim={0 0 0 0.8cm},clip,scale = 0.5]{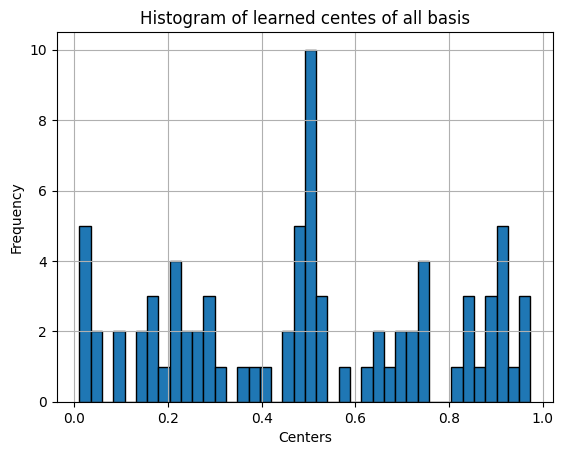}
    \caption{Demonstration of the distributions of the learned basis centers for the Keller-Segel equation. From left to right, top to bottom, we consider models with 20, 40, 60, and 80 basis functions, respectively. The output function exhibits a prominent bump near $x = 0.5$ (see Figure \ref{fig_aggre_ks_sol_demon}), and we observe that more basis functions are concentrated around this region, indicating that the learned basis adapts to the underlying solution structure.}
    \label{fig_ks_center_distribution}
\end{figure}

\textbf{Adaptivity.}
As shown in Figure \ref{fig_aggre_diff_sol_demon}, the terminal solutions exhibit a bump profile and decay rapidly outside the center of the domain.
Therefore, we expect the learned hat functions may focus more closely to the bump.
This behavior is confirmed in our numerical experiments. To validate this hypothesis, we vary the number of basis functions from 40 to 80 and plot the distribution of the learned FEM basis centers. 
See Figure \ref{fig_ks_center_distribution} for the results.

\begin{figure}[H]
    \centering
    \includegraphics[scale = 0.5]{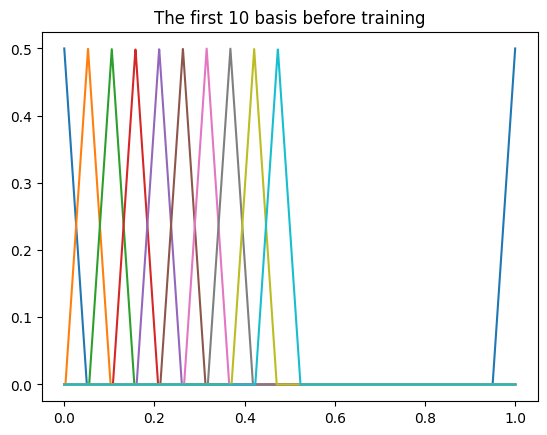}
    \includegraphics[scale = 0.5]{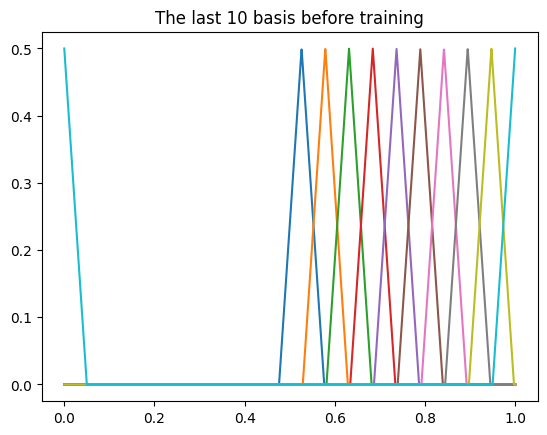}
    \includegraphics[scale = 0.5]{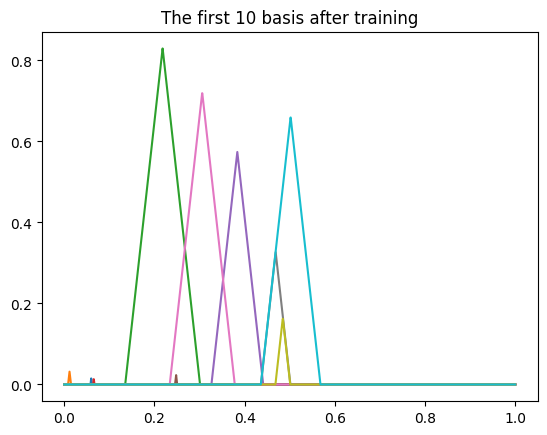}
    \includegraphics[scale = 0.5]{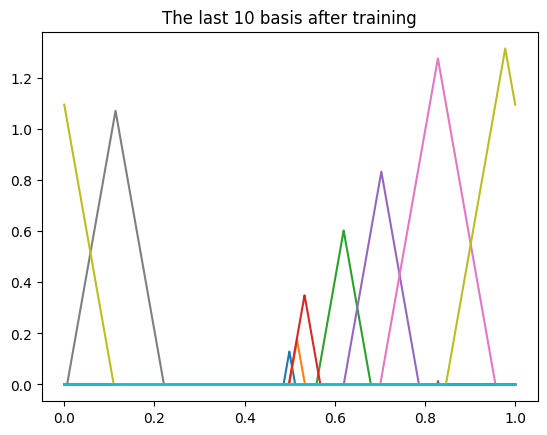}
    \caption{Keller-Segel equation. Demonstration of the FEM basis functions before and after training with 20 basis functions. There is a bump close to $x= 0 .5$ (see Figure \ref{fig_aggre_ks_sol_demon}). 
 We observe from this figure that more localized bases with smaller supports are concentrated around $x = 0.5$ to capture the bump, while fewer bases with larger supports are distributed in other regions to represent the flatter tails (left and right) of the solution.
 }
    \label{fig_ks_basis}
\end{figure}

To better illustrate the concept, Figure \ref{fig_ks_basis} displays the basis functions (20 in total) both before and after training. 
Initially, all basis centers are uniformly distributed with a fixed width $h = 0.05$. 
After training, more basis functions are concentrated near the bump of the terminal solution. The corresponding distribution of basis centers is shown in Figure \ref{fig_ks_center_distribution}.

\subsection{KdV equation} 
\label{sec_kdv}
\begin{figure}[H]
    \centering
    \includegraphics[scale = 0.5]{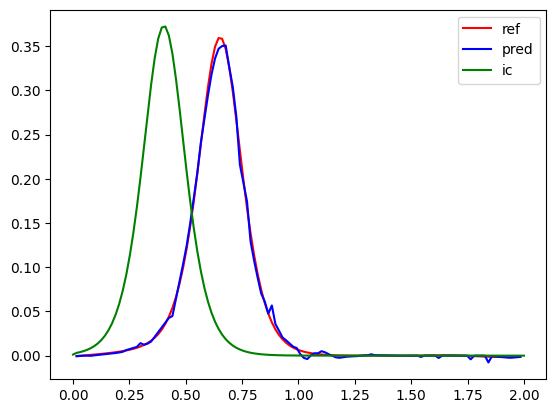}
    \includegraphics[scale = 0.5]{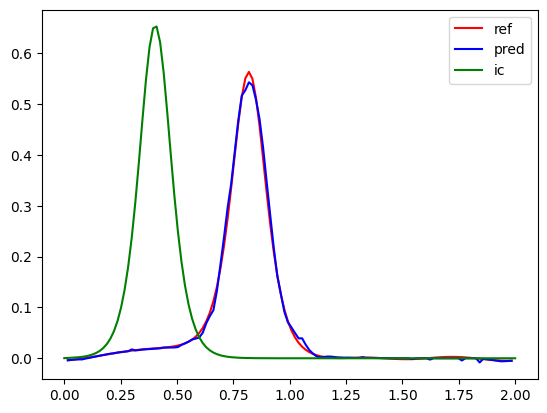}
    \caption{Demonstration of KdV solutions and the prediction using the FEM method with 80 basis functions. The relative errors are $3.86\%$ for the solution shown on the left and $4.23\%$ for the solution on the right.
 }
    \label{fig_kdv_sol}
\end{figure}

In this example, we consider the KdV equation,
\begin{align*}
    u_t - \epsilon u_{xxx} + (u^2/2)_x = 0, x\in[0, 2], t\in [0, 2].
\end{align*}
The target operator maps the initial condition to the solution at the terminal time, with the initial conditions generated using  
$$
u_0(x) = 3c_1\mathrm{sech}(k1(x-x1))^2+3c_2 \mathrm{sech}(k2(x-x2))^2,
$$
with the one free parameter uniformly sampled from $[0, 1]$.
We plot two typical initial conditions and solution pairs in Figure \ref{fig_kdv_sol}.
Notably, all solutions exhibit a bump within the domain $[0, 1]$, and we use this example to further examine the adaptivity of the proposed method.

\begin{figure}[H]
    \centering
    \includegraphics[trim={0 0 0 0.8cm},clip,scale = 0.5]{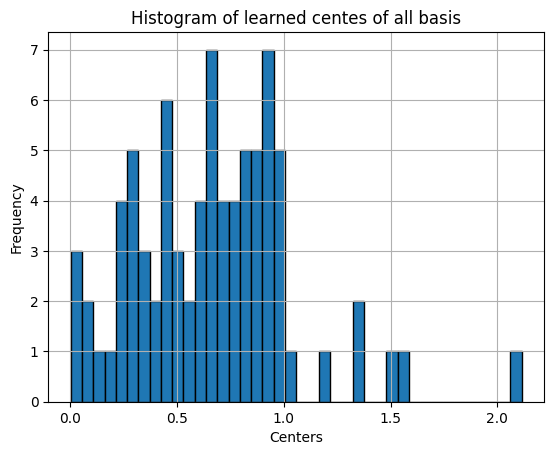}
    \caption{Learned basis center distributions.
    We create a bump in the solution landscape within the region $[0, 1]$ (see Figure \ref{fig_kdv_sol} for examples). As shown in this figure, a greater number of basis functions are relocated to this region (though initialized uniformly distributed in the domain), demonstrating the adaptivity of the proposed method.
 }
    \label{fig_kdv_center}
\end{figure}

\textbf{Setting details.} 
We generate a total of 250 input functions (ICs) and their corresponding solutions for training, where each input function is discretized using 100 sensors (i.e., a mesh of size 100). 
The output functions are evaluated on a mesh with $64$ points for all samples.
All FEM-based models are trained for 2,000 epochs using a cosine annealing learning rate schedule, while DeepONet models are trained with significantly more epochs to reach stable performance. 
We summarize the results in Table \ref{tab_kdv_results}, and present the distribution of the learned centers in Figure \ref{fig_kdv_center}.

In all settings, the coefficient network architecture (used in both DeepONet and FERN) follows a fully connected architecture with dimensions $22 \times 20 \rightarrow 20 \times 1$, and use the hyperbolic tangent ($\tanh$) activation for all layers. 
The centers of the FEM basis functions are uniformly initialized across the domain with a fixed support width of $h = 0.05$. 
For the standard DeepONet, the basis network adopts a deep fully connected architecture with the structure $1 \times 100 \rightarrow 100 \times 100 \rightarrow 100 \times 100 \rightarrow 100 \times 100 \rightarrow 100 \times 100 \rightarrow 100 \times 100 \rightarrow 100 \times N$ with ReLU activations following each hidden layer and a linear output layer, where $N$ is the number of basis; all layers include bias terms except the final one.

\textbf{Accuracy.} As shown in Table \ref{tab_kdv_results}, with the same number of basis functions, compared to DeepONet, FERN gives a smaller error with lower variance and fewer parameters. We also test a smaller DeepONet with a two-layer trunk network while keeping the branch network the same as in all other models. 
Although this reduces the total number of parameters to $113,800$, the relative error increases significantly to $13.25\%\pm4.75\%$.

\textbf{Adaptivity.} As shown in Figure \ref{fig_kdv_sol}, the terminal solutions exhibit a bump in the left half of the domain. We present in Figure \ref{fig_kdv_center} the distribution of learned centers. We observe that most centers concentrate on the left half of the domain, demonstrating the adaptivity of the proposed method.

\begin{table}[H]
\centering
\begin{tabular}{|c|c|c|}
\hline
 & DeepONet(80) & FERN (80) \\
\hline
Relative $L^2$ error         &  $6.0\% \pm 2.81\%$ & $3.93\% \pm 1.10\%$ \\
\hline
$\#$ Coefficient Parameters & 105,600 & 105,600   \\
\hline
$\#$ Basis Parameters        & 58,700  & 160  \\
\hline
$\#$ Total Parameters       & 164,300 & 105,760  \\
\hline 
\end{tabular}
\caption{KdV results with $80$ basis.
We also test a smaller DeepONet with a two-layer trunk network while keeping the branch network the same as in all other models. 
Although this reduces the total number of parameters to $113,800$, the relative error increases significantly to $13.25\%\pm4.75\%$.
}
\label{tab_kdv_results}
\end{table}

\subsection{Viscous Burgers' equation}

\begin{figure}[H]
    \centering
    \includegraphics[scale = 0.5]{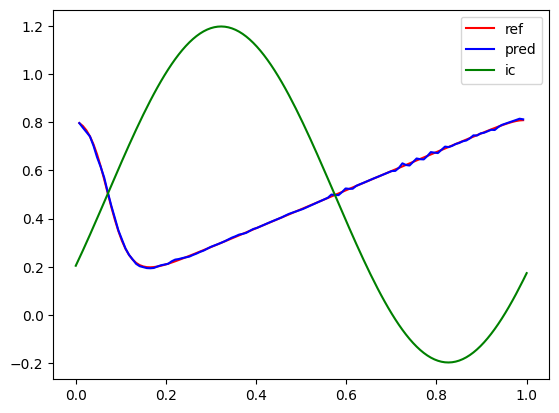}
    \includegraphics[scale = 0.5]{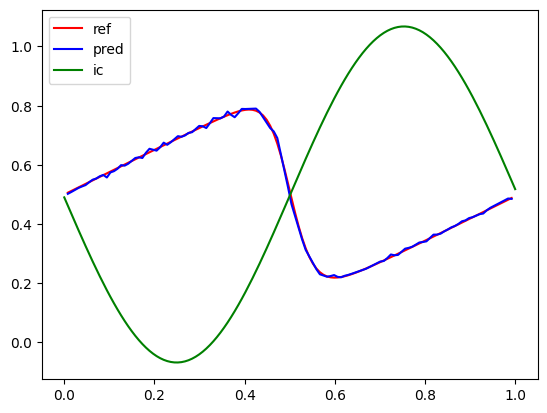}
    \caption{Demonstration of viscous Burgers' equation solutions and the prediction using the FEM method with 40 basis functions. The relative errors are $0.71\%$ for the solution shown on the left and $1.01\%$ for the solution on the right.
 }
    \label{fig_vburgers_sol}
\end{figure}

Lastly, as a widely tested example in operator learning, we consider the Viscous Burgers' equation,
\begin{align*}
    u_t - 0.01 u_{xx} + (u^2/2)_x = 0, x\in[0, 1], t\in [0, 1].
\end{align*}
The target operator maps the initial condition to the solution at the terminal time, with the initial conditions generated using 
$$
u_0(x) = c_1\sin(2\pi(x-c_0)) + 0.5,
$$
with 
$c_0$ uniformly sampled from $[0, 0.5]$ and 
$c_1$ uniformly sampled from $[0.5, 1.0]$.
We plot two typical initial conditions and solutions pairs in Figure \ref{fig_vburgers_sol}.

\textbf{Setting details.} 
We generate a total of 84 input functions (ICs) and their corresponding solutions for training, where each input function is discretized using 22 sensors (i.e., a mesh of size 22). 
The output functions are evaluated on a mesh with $64$ points for all samples.
All FEM-based models are trained for 2,000 epochs using a cosine annealing learning rate schedule, while DeepONet models are trained with significantly more epochs to reach stable performance. 
We summarize the results in Table \ref{tab_vburgers_results}.

In all settings, the coefficient network architecture (used in both DeepONet and FERN) follows a fully connected architecture with dimensions $22 \times 20 \rightarrow 20 \times 1$, and use the hyperbolic tangent ($\tanh$) activation for all layers. 
The centers of the FEM basis functions are uniformly initialized across the domain with a fixed support width of $h = 0.05$. 
For the standard DeepONet, the basis network adopts a deep fully connected architecture with the structure $1 \times 100 \rightarrow 100 \times 100 \rightarrow 100 \times 100 \rightarrow 100 \times 100 \rightarrow 100 \times 100 \rightarrow 100 \times 100 \rightarrow 100 \times K$, with ReLU activations following each hidden layer and a linear output layer; all layers include bias terms except the final one.

\textbf{Results analysis.} As shown in Table \ref{tab_vburgers_results}, the propsoed method gives an accuracy comparable to DeepONet but with fewer parameters.

\begin{table}[H]
\centering
\begin{tabular}{|c|c|c|}
\hline
 & DeepONet(40) & FERN (40) \\
\hline
Relative $L^2$ error         &  $0.95\% \pm 0.41\%$ & $0.93\% \pm 0.31\%$ \\
\hline
$\#$ Coefficient Parameters & 19,200 & 19,200   \\
\hline
$\#$ Basis Parameters        & 54,700  & 80  \\
\hline
$\#$ Total Parameters       & 73,900 & 19,280  \\
\hline 
\end{tabular}
\caption{Viscous Burgers example results with $40$ basis.
}
\label{tab_vburgers_results}
\end{table}

\section{Discussion}
\label{sec_conclusion}
In this paper, we developed a finite-element polynomial basis–based operator learning framework that efficiently solves families of PDEs with a substantially reduced number of trainable parameters. By constructing adaptive FEM bases through a shallow neural architecture, the proposed method retains the flexibility of learnable representations while inheriting the local adaptivity and interpretability of traditional finite element methods. Numerical experiments across seven distinct PDE families confirm that the method achieves competitive or superior accuracy compared to standard DeepONet models, often with orders of magnitude fewer parameters and smaller variance in prediction errors. 
The proposed framework opens several promising directions for future research. One avenue is to construct higher-order FEM bases by incorporating different polynomial families (e.g., Legendre, Chebyshev, or hierarchical bases) and exploring alternative activation functions to enhance expressiveness. Another important direction is to analyze and improve extrapolation capabilities, as the adaptive nature of the learned FEM bases suggests potential for stronger generalization beyond the training regime. Overall, this study bridges classical finite element analysis and modern operator learning, offering a pathway toward efficient, interpretable, and adaptive neural operator frameworks for solving complex PDE systems.

\section*{Acknowledgment}
ZZ was supported by the U.S. Department of Energy (DOE) Office of Science Advanced Scientific Computing Research program DE-SC0025440.
HL was partially supported by HKRGC ECS grant 22302123, HKRGC GRF 12301925, and Guangdong and Hong Kong Universities ``1+1+1” Joint Research Collaboration Scheme UICR0800008-24.
HS was supported in part by NSF DMS 2427558. 
GL was supported by the NSF (DMS-2053746, DMS-2134209, ECCS-2328241, CBET-2347401 and OAC-2311848), and U.S.~Department of Energy (DOE) Office of Science Advanced Scientific Computing Research program DE-SC0023161, the Uncertainty Quantification for Multifidelity Operator Learning (MOLUcQ) project (Project No. 81739), and DOE–Fusion Energy Science, under grant number: DE-SC0024583.

\bibliographystyle{siam}
\bibliography{references}
\end{document}